\documentclass[a4paper]{article}

\usepackage[a4paper]{geometry}
\usepackage{amsmath} 
\usepackage{amssymb}
\usepackage{latexsym}
\usepackage{dsfont}
\usepackage{euler}
\usepackage{enumerate} 
\usepackage{xspace} 
\usepackage{url} 
\usepackage{natbib}
\usepackage{booktabs} 

\usepackage{ntheorem} 
\theorembodyfont{\upshape}

\usepackage[open-square]{QED} 
\usepackage{longtable} 
\usepackage{sectsty} 

\usepackage{tikz} 
\usetikzlibrary{arrows,decorations.markings}
\tikzset{>=triangle 45}

\usepackage{pf2} 
\pflongnumbers 

\usepackage[affil-it]{authblk} 

\newcommand{\anondivision}{\begin{center}*\hspace{2em}*\hspace{2em}*\end{center}}
\renewcommand{\iff}{\Leftrightarrow}
\newcommand{\imp}{\Rightarrow}
\newcommand{\of}[1]{\left({#1}\right)}

\newcommand{\setof}[1]{\left\{{#1}\right\}}

\newcommand{\abs}[1]{\left|{#1}\right|}
\newcommand{\norm}[1]{\left\|{#1}\right\|}
\newcommand{\pres}[2]{\left\langle #1 \middle| #2 \right\rangle}

\newcommand{\D}[1]{\mathrm{d}{#1}}
\newcommand{\choice}[1]{\left\{\begin{array}{cl}#1\end{array}\right.}
\newcommand{\supp}{\text{supp}}
\newcommand{\Sc}[1]{\mathscr{#1}}
\newcommand{\dotcup}{\,\ensuremath{\mathaccent\cdot\cup}\,}
\newcommand{\fn}[2]{#1\of{#2}}
\newcommand{\defeq}{:=}
\newcommand{\Folner}{F\o lner\xspace}

\newtheorem{theorem}{Theorem}[section]
\newtheorem{proposition}[theorem]{Proposition}
\newtheorem{lemma}[theorem]{Lemma}
\newtheorem{corollary}[theorem]{Corollary}
\newtheorem{theorem-ac}[theorem]{Theorem (AC)}
\newtheorem{proposition-ac}[theorem]{Proposition (AC)}
\newtheorem{lemma-ac}[theorem]{Lemma (AC)}
\newtheorem{corollary-ac}[theorem]{Corollary (AC)}
\newtheorem{remark}[theorem]{Remark}
\newtheorem{definition}[theorem]{Definition}
\newtheorem{example}[theorem]{Example}

\newtheorem{question}[theorem]{Question}
\newtheorem{definition-theorem}[theorem]{Definition/Theorem}

\title{\bfseries\sffamily Fair amenability for semigroups}
\author{Josh Deprez\thanks{Electronic address: \texttt{josh.deprez@gmail.com}.}}
\affil{School of Mathematics \& Physics \\ University of Tasmania}

\begin{document}
\maketitle

\begin{abstract}
A new flavour of amenability for discrete semigroups is proposed that generalises group amenability and follows from a \Folner-type condition. Some examples are explored, to argue that this new notion better captures some essential ideas of amenability. A semigroup $S$ is left fairly amenable if, and only if, it supports a mean $m\in\ell^\infty(S)^*$ satisfying $\fn{m}{f} = \fn{m}{s\ast f}$ whenever $s\ast f\in\ell^\infty(S)$, thus justifying the nomenclature ``fairly amenable''.
\end{abstract}


\section{Introduction}

Amenability begin in essence alongside modern analysis, as it is a central property lacking in a group used to show, for example, the Banach-Tarski paradox \citep{Wagon:1993th}. The first working definition for what is now called amenability was given by \cite{VonNeumann1929}, in terms of finitely-additive measures. A group $G$ is amenable if there is a finitely-additive measure $\mu$ such that $\fn\mu{G} = 1$, and $\fn\mu{gA} = \fn\mu{A}$ for all $g\in G, A\subseteq G$ ($\mu$ is \emph{left invariant}). This definition has the advantages of being easy to comprehend, hiding very little, and it is easy to show that the free group on two generators $\mathbb{F}_2$ does not support such a finitely-additive measure.

The first modern definition of amenability, in its form as extended to semigroups, was given by \cite{Day1956}, whose concept involved invariant means. A \emph{mean} is a non-negative linear functional $m\in\ell^\infty(S)^*$ such that $\fn{m}{\chi_S} = 1$. The means generalise the finitely-additive measures: to obtain a mean from a finitely-additive measure, use the Lebesgue integral construction. An element $s\in S$ acts on a function $f\in\ell^\infty(S)$ (on the left), by setting $\fn{(s\cdot f)}{t} \defeq \fn{f}{st}$ for all $t\in S$. Briefly, then, a semigroup $S$ is (classically) \emph{left amenable} when there exists such an $m$ satisfying
\(
    \fn{m}{s\cdot f} = \fn{m}{f}
\)
for all $f\in\ell^\infty(S)$. For groups this coincides exactly with the von Neumann condition. However, for semigroups it does not: in fact, a mean satisfies the above condition if, and only if, the associated finitely-additive measure $\mu$ satisfies\footnote{Recall that for any $s\in S$ and $A\subseteq S$, the set $s^{-1}A \defeq \setof{t\in S: st \in A}$.} $\fn\mu{s^{-1}A} = \fn\mu{A}$ for all $s\in S$ \citep{AlanL.T.Paterson1988}. This might be called left \emph{preimage} invariance of $\mu$. A simple but surprising consequence of all this is that all semigroups with a zero element are both left and right amenable \citep{Day1956} yet they cannot have a  (totally) invariant finitely-additive measure \citep[p231]{vanDouwen:1992ux}. On the other hand, all semigroups with more than one distinct left zero are not left amenable \citep{AlanL.T.Paterson1988}.

Numerous other alternative definitions for amenability from group theory disagree on semigroups in general. The \Folner conditions, originally shown for groups by \cite{Folner1955} and of which there are now several flavours, have varying degrees of relation to left amenability of a semigroup. The \Folner conditions are useful for showing when a group has amenability, and effectively describe the essential reason all Abelian groups are amenable. \Folner's original conditions were first generalised to semigroups by Frey in 1960 and subsequently a simpler proof was given by \cite{Namioka1964}. Some of the \Folner-type criteria that are sufficient for left amenabilty of a semigroup include the weak and strong \Folner conditions \citep{Argabright:1967vm} and the weak and strong \Folner-Namioka conditions \citep{Yang:1987ta}. A \emph{necessary} \Folner-type condition for amenable semigroups is the one described by \cite{Namioka1964}. 

For some of these \Folner conditions, and other related conditions, if the semigroup in question is \emph{cancellative}, then there are improved results, since the inequality $2\abs{A\backslash sA} \ge \abs{sA\,\triangle\,A} \ge 2\abs{sA\backslash A}$, true for any $s\in S$ and finite $A\subseteq S$, is then saturated. For example, Frey's thesis showed that if $S$ is a cancellative semigroup that contains no free subsemigroup on two generators, and is left amenable, then every subsemigroup of $S$ is left amenable. An improvement was made recently by \cite{Donnelly2012}: if $T$ is a subsemigroup of $S$, $S$ is cancellative, $T$ does not contain a free subsemigroup on two generators, and $S$ is left amenable, then $T$ is left amenable. 

Another set of results concerns translating amenability between groups and algebras. A Banach algebra $\mathfrak{A}$ is called \emph{amenable} if $\mathcal{H}^1(\mathfrak{A},E^*) = \setof{0}$ for every Banach $\mathfrak{A}$-bimodule $E$ \citep[p43]{Runde:th}---this is equivalent to saying all derivations are inner derivations. It is the famous theorem of \cite{Johnson:1972ui} that shows that the group $G$ is amenable if, and only if, $\ell^1(G)$ is amenable (as a convolution Banach algebra). However, for a semigroup $S$, the amenability of $\ell^1(S)$ does not relate well to the amenability of $S$.

\subsection{Inverse semigroups}

One might hope that the situation would be less complicated when restricted to the class of inverse semigroups. Sticking to classical amenability, it is so much less complicated as to be almost trivial: \cite{Duncan1978} showed that an inverse semigroup $S$ is amenable if and only if its maximal group homomorphic image (denoted $G(S)$) is amenable. As an example, if the inverse semigroup $S$ has a zero, then $G(S)$ is the trivial group, and therefore $S$ is amenable. 

On the other hand, the convolution Banach algebra $\ell^1(S)$ is amenable if, and only if, the semilattice of idempotents (denoted $E(S)$) is finite and every subgroup of $S$ is amenable. This is regarded as too strong \citep{Milan}, since it eliminates many commutative inverse semigroups.

\cite{Paterson1998} suggested the following result points at one resolution: if the inverse semigroup $S$ has all maximal subgroups amenable, then $VN(S)$ (the von Neumann algebra of $S$) is amenable. 

\cite{Milan:2008tp} argued that the weak containment property---another generalisation of amenability for groups---is an appropriate notion of amenability for inverse semigroups, by showing the following. The free group on two generators with a zero adjoined, an example of a Clifford semigroup, does not have weak containment, but the commutative inverse semigroups all have weak containment. Therefore the weak containment property sits neatly between amenability of $S$ and amenability of $\ell^1(S)$. \cite{Milan:2008tp} also showed that an $E$-unitary inverse semigroup has weak containment if, and only if, $G(S)$ is amenable, and that examples of inverse semigroups with weak containment include the graph inverse semigroups, which generalise and include the polycyclic monoids \citep[see][]{Jones2011}. 

Recall that for any given inverse semigroup $S$, the left regular representation $\pi_2$ of $s\in S$ on the Hilbert space $\mathcal{H} = \ell^2(S)$ is defined by
\[
    \fn{\pi_2}{s}f \defeq \sum_{tt^* \le s^*s} \fn{f}{t} st \quad\text{for all }f \in \ell^2(S)
\]
\citep{Paterson1998}. This representation is central to the weak containment property. Due to the reliance the natural partial order to keep the above summation well-defined (consider $\fn{\pi_2}{0}f$: the only idempotent bounded above by $0$ is $0$), this may not be adequately generalisable to arbitrary semigroups.

\anondivision

In the remainder of this paper I describe a condition, similar to amenability, and given in terms of finitely-additive measures, which was inspired by the results relating to cancellative semigroups and the regular representations of an inverse semigroup, that takes advantage of zeroes and other non-cancellative elements in a natural way. An extension of this condition to the context of means will be given in a forthcoming paper. The Axiom of Choice shall be assumed throughout, though it will be mentioned where used.


\section{Definitions}


The following is required to introduce the condition. Let $S$ be a semigroup, and define the maps
\[
    \lambda_s(x) \defeq sx;\quad \rho_s(x) \defeq xs\quad \text{for all }s,x\in S.
\]    
$\lambda$ and $\rho$ are known as the \emph{left regular} and \emph{right regular} representations, respectively. (Note that these should not be confused with the regular representations on a Hilbert space described above.) For all $s\in S$, $\lambda_s$ and $\rho_s$ are elements of $\mathcal{T}_S$, the transformation semigroup of the set $S$.

\begin{definition}[Acting injectively]\label{def:injectiveaction}
If $\lambda_s|_A:A\to sA$ is an injection, then $s$ is said to \emph{act injectively on the left of $A$}. If $\rho_s|_A: A\to As$ is an injection, then $s$ acts injectively on the \emph{right} of $A$.
\end{definition}

By definition, every $s\in S$ acts injectively on the left of $S$ if, and only if, $S$ is left cancellative. Finding the sets acted on injectively permits analysing any semigroup, rather than unsubtly requiring the semigroup to be cancellative. 

\begin{lemma}
For any $s\in S$ and $A\subseteq S$, the following are equivalent:
\begin{enumerate}[(i)]
\item $s$ acts injectively on the left of $A$;
\item For all two-element set $F\subseteq A$, $\abs{sF} = \abs{F}$;
\item For any finite set $F\subseteq S$, $\abs{s(F\cap A)} = \abs{F\cap A}$.
\end{enumerate}
\end{lemma}


\begin{definition}[Subinvariant]
Let $S$ be a semigroup, and $\mu$ a finitely-additive measure on $S$ with finite total measure. If  
\[
	\fn\mu{sA} \le \fn\mu{A} \quad [\fn\mu{As} \le \fn\mu{A}] \quad \text{for all }s\in S\text{ and } A\subseteq S,
\]
then we say $\mu$ is \emph{left [right] sub-invariant}.
\end{definition}

Suppose $sA = A$ (for instance $s$ is an identity); it is then clear that the inequality above cannot be strict in general. 

Suppose for some element $s$ and set $A$ there is some $s'$ such that $s'sA = A$. Then both $s$ and $s'$ are behaving as injective acts, and when restricted to $A$, $s's$ acts as a permutation of $A$. If $\mu$ is left sub-invariant, 
\[
	\fn\mu{A} = \fn\mu{s'sA} \le \fn\mu{sA} \le \fn\mu{A},
\]
and thus $\fn\mu{sA} = \fn\mu{A}$. This suggests the next definition, which is the most important here.

\begin{definition}[Fairly invariant, fairly amenable]\label{def:fairly amenable}
Let $S$ be any semigroup, let $\mu$ a finitely-additive measure on $S$ with $\fn\mu{S} = 1$, and let $s\in S$ and $A\subseteq S$. 
 
If whenever $s$ acts injectively on the left [right] of $A$,
\[
	\fn\mu{sA} = \fn\mu{A} \quad [\fn\mu{As} = \fn\mu{A}]
\]
then $\mu$ is \emph{fairly left [right] invariant}. If such a $\mu$ exists for a given semigroup $S$, then $S$ is \emph{fairly left [right] amenable}.
\end{definition}
In other words, invariance of $\mu$ is only required in the places where an element $s$ acts injectively on the set. As we shall see, this weakening of total invariance handles the issue discussed in \cite[p231]{vanDouwen:1992ux}.

\begin{figure}
\begin{center}
\begin{tikzpicture}[scale=0.7]
    \filldraw[draw=red,fill=red!10] (0,0) circle (3) node[anchor=north,below=1.3cm] {$A$};
    \filldraw[draw=blue,fill=blue!10] (0,0) circle (1.5) node {$B$};
    \filldraw[draw=violet,fill=violet!10] (6,0) circle (1.5) node {$sA = sB$};
    \filldraw[fill=black] (6,0.7) circle (0.07) node[anchor=west] {$x$}
        (0,0.7) circle (0.07) node[anchor=south] {$b$}
        (0.5,2) circle (0.07) node[anchor=south] {$a$}
        (3,2.8) circle (0.07) node[anchor=south] {$c\in s^{-1}\setof{x}$};
    \begin{scope}[decoration={markings,mark=at position 0.5 with {\arrow{>}}}]
        \draw[red,postaction={decorate}] (0.7,2.92) -- node[anchor=north]{$s$} (6.4,1.45);
        \draw[red,postaction={decorate}] (0.7,-2.92) -- node[anchor=south]{$s$} (6.4,-1.45);
        \draw[blue,postaction={decorate}] (0,1.5) -- node[anchor=south]{$s$} (6,1.5);
        \draw[blue,postaction={decorate}] (0,-1.5) -- node[anchor=south]{$s$} (6,-1.5);
        \draw[postaction={decorate}] (0, 0.7) -- (6,0.7);
        \draw[postaction={decorate}] (0.5, 2) -- (6,0.7);
        \draw[postaction={decorate}] (3,2.8) -- (6,0.7);
    \end{scope}
\end{tikzpicture}
\caption{\label{fig:theobviouslemma} For every set $A$ and element $s$ there is some subset $B$ such that $sA = sB$ and $s$ acts injectively on $B$.}
\end{center}
\end{figure}
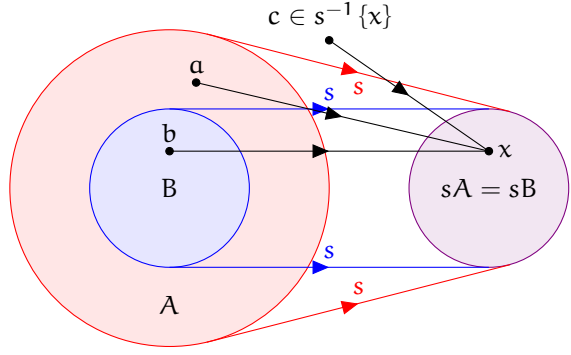

\begin{lemma} For any semigroup $S$ and finitely-additive probability measure $\mu$, left [right] fair invariance of $\mu$ implies left [right] sub-invariance of $\mu$.
\Proof{For a pictoral overview see Figure \ref{fig:theobviouslemma}.
\begin{pruf}
\step{1}{For any $A\subseteq S$ and $s\in S$ there exists a $B\subseteq A$ such that $sA = sB$ and $s$ is injective on $B$.}
    \begin{pruf}
    \pf\ Use the Axiom of Choice to choose one $b\in s^{-1}\setof{x}\cap A$ for each $x \in sA$. $B$ is simply the set of those choices.
    \end{pruf}
\step{2}{If $B\subseteq A\subseteq S$, and $sA = sB$, and $s$ acts injectively on $B$ (but not necessarily on $A$), then $\fn\mu{A} \ge \fn\mu{sA}$.}
    \begin{pruf}
    \pf\begin{align*}
            \fn\mu{A} &\ge \fn\mu{B} \quad\because B\subseteq A\\
                &= \fn\mu{sB} \quad\because\text{fair invariance of }\mu\\
                &= \fn\mu{sA} \quad\because sB = sA.
        \end{align*}
        as required.\qed
    \end{pruf}
\end{pruf} }
\end{lemma}

\begin{remark}
What about selecting $\fn\mu{sA} \ge \fn\mu{A}$ as a condition (``super-invariance'')? If $sA$ is a subset of $A$ then $\fn\mu{sA} = \fn\mu{A}$, and so disjoint subsets $sA, tA$ may lead to a contradiction.

By definition, if $s$ acts injectively on the left of $A$, then $sa = sb \imp a = b$ for any $a,b\in A$ and $s$ is \emph{left cancellative on $A$}; similarly on the right. Hence another way of defining fair invariance is in terms of cancellation. Groups are totally cancellative both ways, but there are non-group examples of left- and right-cancellative semigroups.
\end{remark}


\section{Consequences}

\begin{figure}
\begin{center}
\begin{tikzpicture}
    \begin{scope}[fill opacity=0.2,text opacity=1]
        \filldraw[fill=red,draw=black] (0,0) node[anchor=north] {$A$} circle (1);
        \filldraw[fill=red,draw= black] (3,0) node[anchor=north] {$sA$} circle (1);
        \filldraw[fill=blue,draw= black] (0.5,1.3) node[anchor=south] {$F_n$} circle (1);
        \filldraw[fill=blue,draw= black] (3.5,1.0) node[anchor=south] {$sF_n$} circle (0.6);
    \end{scope}
    \begin{scope}[decoration={markings,mark=at position 0.5 with {\arrow{>}}}]
        \draw[postaction={decorate}] (0,-1) -- node[anchor=south] {$s$} (3,-1);
        \draw[postaction={decorate}] (0,1) -- (3,1);
        \draw[postaction={decorate}] (0.5,0.3) -- (3.5,0.4);
        \draw[postaction={decorate}] (0.66,2.29) -- node[anchor=north] {$s$} (3.65,1.58);
    \end{scope}
\end{tikzpicture}
\caption{\label{fig:injectiveonintersection}If $s$ acts injectively on $A$, then it also acts injectively on the subset $A\cap F_n$ of $A$, and so $\abs{A\cap F_n} = \abs{s(A\cap F_n)}$. Note that $s(A\cap F_n)\subseteq sA \cap sF_n$ might not be saturated---consider disjoint $A$ and $F_n$.}
\end{center}
\end{figure}
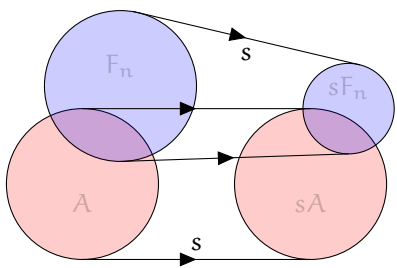

Fair amenability is a generalisation of amenability for groups, as follows.
\begin{corollary}\label{cor:backwardscompatible}
A group is amenable if, and only if, it is fairly amenable. 
\Proof{This is trivial since every element $g$ in a group $G$ acts bijectively on $G$, and so a finitely-additive measure on $G$ is invariant if, and only if, it is fairly invariant. \qed}
\end{corollary}
Similar to classical amenability, fair amenability is also a consequence of a \Folner-type condition, as follows.
\begin{theorem}\label{theorem:thenewfolner}
Let $S$ is a countable semigroup. If for each $s\in S$ there exists a sequence of non-empty finite sets $\setof{F_n}_{n\in\mathbb{N}}$ eventually covering $S$ such that for all $A\subseteq S$,
\[
    \lim_{n\to\infty} \frac{\abs{s(A\cap F_n)\,\triangle\, (sS \cap F_n)} }{\abs{F_n}} = 0,
\]
then $S$ is left fairly amenable. (Similarly for $F_ns$/on the right.)
\Proof{
Fix a free ultrafilter $U$ over $\mathbb{N}$ and define $\mu$ through the ultralimit
\[
    \fn\mu{A} \defeq \lim_U \frac{\abs{A\cap F_n}}{\abs{F_n}} \quad\text{for all } A\subseteq S.
\]
\begin{pruf}
\step{fapm}{For any set $A$, the ultralimit above exists, and $\mu$ is a finitely-additive measure with $\fn\mu{S} = 1$.}
    \begin{pruf}
    \pf\ 
    The sequence is bounded so by the Bolzano-Weierstra\ss{} Theorem there is always a convergent subsequence, so the ultralimit always exists. $\mu(S) = 1$ since $\abs{S\cap F_n} = \abs{F_n}$ for all $n$, in which case the sequence is constantly $1$. $\mu$ is finitely-additive as a simple consequence of \L os's Theorem, in particular, $\lim_U (x_n + y_n) = \lim_U x_n + \lim_U y_n$ for any sequences $\setof{x_n}_{n\in\mathbb{N}}, \setof{y_n}_{n\in\mathbb{N}}$ where the ultralimits exist.
    \end{pruf}
\step{invariance}{$\mu$ is left fairly invariant.}
    \begin{pruf}
    \pf\ Suppose $s$ acts injectively on the left of $A$, and thus every subset of $A$, in particular $A\cap F_n$. Then $\abs{A\cap F_n} = \abs{s(A\cap F_n)}$. See Figure \ref{fig:injectiveonintersection}. Then,
    \begin{align*}
        \abs{ \frac{\abs{A\cap F_n}}{\abs{F_n}} -  \frac{\abs{sA\cap F_n}}{\abs{F_n}} } &= \frac{\abs{ \abs{A\cap F_n} -  \abs{sA\cap F_n} }}{\abs{F_n}} \\
            &= \frac{\abs{ \abs{s(A\cap F_n)} - \abs{sA\cap F_n} }}{\abs{F_n}} \\
            &\le \frac{\abs{s(A\cap F_n)\,\triangle\, (sS \cap F_n)}}{\abs{F_n}} \to 0 \text{ as } n\to\infty
    \end{align*}
    by hypothesis, and hence $\fn\mu{A} = \fn\mu{sA}$, as required. \qed
    \end{pruf}
\end{pruf} }
\end{theorem}

\begin{remark}
While there are semigroups lacking strong \Folner sequences that are also fairly amenable, this appears to be mitigated in this condition as the top line is contained in the right ideal $sS$. Consider, for example, an infinite amenable group $G$ with zero adjoined ($G^0$), which is fairly amenable (see Corollary \ref{cor:addingazero} below) and the zero element has no associated \Folner sequence, however $0S = \setof{0}$ and therefore any \Folner sequence will do in Theorem \ref{theorem:thenewfolner}.
\end{remark}

\begin{corollary}\label{cor:allfinite}
All finite semigroups $S$ are fairly amenable (both ways). 
\Proof{
$s$ is injective on the left of $A \subseteq S$ if, and only if, $\abs{sA} = \abs{A}$; similarly on the right. Therefore the counting measure suffices. Alternatively, use the constant \Folner sequence $\setof{S}_{n\in\mathbb{N}}$: for any $A\subseteq S$, 
	\[
		\frac{\abs{s(A\cap S)\,\triangle\,(sA\cap S)}}{\abs{S}} = \frac{\abs{sA\,\triangle\,sA}}{\abs{S}} = \frac{0}{\abs{S}} = 0
	\]
    as required by Theorem \ref{theorem:thenewfolner}. \qed}
\end{corollary}

\begin{remark}
Suppose that, given some set $A$, $\fn\mu{sA} = \fn\mu{A} [\fn\mu{As} = \fn\mu{A}]$ for any $s$. We may describe $A$ as being a left [right] $\mu$-invariant set. In a fairly left [right] amenable semigroup $S$, every singleton set $\setof{x}$ for $x\in S$ is guaranteed to be a left [right] invariant set.
\end{remark}

\begin{lemma}\label{lem:finitemassless}
Let $S$ be an infinite left [right] fairly amenable semigroup with measure $\mu$, having a left [right] zero $z\in S$. If $F$ is a finite subset of $S$, then $\fn\mu{F} = 0$.
\Proof{
\begin{pruf}
\step{allsame}
{Every singleton set has the same measure $k$.}
    \begin{pruf}
    \pf\ We can go via $\setof{z}$: for any $s,t\in S$, 
    \[
        \fn\mu{\setof{s}} = \fn\mu{z\setof{s}} = \fn\mu{\setof{z}} = \fn\mu{z\setof{t}} = \fn\mu{\setof{t}}.
    \] 
    \end{pruf}
\step{k0}{$k=0$, therefore $\fn\mu{F} = 0$.}
    \begin{pruf}
    \pf\ If $k>0$ there exists some finite $N$ such that $Nk > 1$, i.e. the disjoint union of $N$ singletons would have greater than 1 measure. Hence $k = 0$. Then
    \[
        \fn\mu{F} = \sum_{f\in F}\fn\mu{\setof{f}} = \sum_{f\in F} k = 0.
    \]
    The right case holds similarly. \qed
     \end{pruf}
\end{pruf}}
\end{lemma}

\begin{corollary}\label{cor:diracmeasure}
Let $S$ be a non-trivial semigroup with zero. The finitely-additive measure $\delta_0$ given by
\[
	\fn{\delta_0}{A} = \fn{\delta_0}{0^{-1}A} = \fn{\delta_0}{A0^{-1}} = \choice{1&\text{if } 0\in A\\ 0&\text{if } 0\notin A}
\]
(i.e. the measure obtained from the invariant mean $m\in\ell^\infty(S)^*$ given by $\fn{m}{f} = \fn{f}{0}$ for all $f\in\ell^\infty(S)$) cannot be fairly invariant.
\Proof{
Let $a\in S$ where $a\ne 0$ and assume $\delta_0$ is fairly invariant. Then
\begin{align*}
	1 &= \fn{\delta_0}{\setof{0}} \quad\text{by definition}\\
	  &= \fn{\delta_0}{0\setof{a}} \\
	  &= \fn{\delta_0}{\setof{a}} \\
	  &= 0 \quad\text{since } 0\notin\setof{a},\\
\end{align*}
contradiction. 
\qed}
\end{corollary}

\begin{question}Is there a left or right fairly amenable infinite semigroup with a finite subset having positive mass?
\end{question}

\anondivision

Recall the \emph{Green's relations} $\Sc{L}, \Sc{R}, \Sc{D}, \Sc{H}, \Sc{J}$ on a semigroup. There are two easy lemmas. 
\begin{lemma}\label{lem:finitemassless2}
If $S$ is left [right] fairly amenable with measure $\mu$, any finite subset $F$ of an infinite $\Sc{L}$-class [$\Sc{R}$-class] has $\fn\mu{F} = 0$. It follows that in either case any finite subset $F$ of a $\Sc{H}$-class has $\fn\mu{F} = 0$, and if $S$ is fairly amenable on both sides than any finite subset $F$ of a $\Sc{D}$-class has $\fn\mu{F} = 0$.
\Proof{
\begin{pruf}
\step{allrelated}{Every singleton subset of an $\Sc{L}$-class has the same measure $k$.}
    \begin{pruf}
    \pf\ By definition, for all $a,b\in S$ such that $a\,\Sc{L}\,b$, there exists $s,s'\in S^1$ such that $sa=b, s'b=a$, and we only need one of these to establish that if $\mu$ is the left fairly invariant finitely-additive measure, 
    \[
        \fn\mu{\setof{a}} = \fn\mu{s\setof{a}} = \fn\mu{\setof{sa}} = \fn\mu{\setof{b}} \quad\text{for all } a,b\in S.
    \] 
    \end{pruf}
\step{allzero}{Every finite subset has measure 0.}
    \begin{pruf}
    \pf\ As for the final step of Lemma \ref{lem:finitemassless}. \qed
    \end{pruf}
\end{pruf}
}
\end{lemma}

Green's Lemma \citep[p43]{Howie1976} states that for any $a,b\in S$ such that $a\,\Sc{R}\,b$, the restricted right regular representations to $\Sc{L}$-classes, $\rho_s|L_a$ and $\rho_{s'}|L_b$, are mutually inverse $\Sc{R}$-class preserving bijections between the $\Sc{L}$-classes $L_a$ and $L_b$. Put another way, there exists an $s\in S$ that acts injectively on the right of $L_a$ and an $s'\in S$ that acts injectively on the right of $L_b$.

\begin{lemma}
Let $S$ be a semigroup. If $S$ is right fairly amenable with measure $\mu$ then within each $\Sc{D}$-class all $\Sc{L}$ classes have the same measure. Similarly, if $S$ is left fairly amenable with $\mu$ then within each $\Sc{D}$-class all $\Sc{R}$ classes have equal measure. It follows that if $S$ is fairly amenable (both ways) then all $\Sc{D}$-related $\Sc{H}$-classes have equal measure.
\Proof{ 
Suppose $L_a, L_b$ are $\Sc{L}$-classes contained within the same $\Sc{D}$-class.
\begin{pruf}
\step{exists}{There exist $s,s' \in S^1$ such that $L_a = L_bs'$ and $L_b = L_as$ are both examples of injective right acts.}
    \begin{pruf}
    \pf\ Use Green's Lemma.
    \end{pruf}
\step{qed}{$\fn\mu{L_a} = \fn\mu{L_as} = \fn\mu{L_b}$. \qed}
\end{pruf}
}
\end{lemma}

What can we say about the value of a fairly invariant finitely-additive measure $\mu$ between distinct $\Sc{D}$-classes? Probably not a lot (see Example \ref{ex:cliffordF2andF1} below). 
\anondivision
A result for groups states that the direct product of finitely many amenable groups is also amenable. This is easily shown by noting that if $G = G_1\times G_2$ then the subgroup $H = \setof{(g_1, 1_{G_2}): g_1\in G_1} \cong G_1$, and $G/H\cong G_2$, so therefore the amenability of $G_1$ and $G_2$ imply the amenability of $H$ and $G/H$, and hence $G$. The fair amenability analogue of this result is as follows, but shown in a more involved manner.

\begin{theorem}\label{thm:directproductfairly}
Let $S, T$ be semigroups that are each left [right] fairly amenable. $S\times T$ is as well.
\Proof{
Let $\mu_S$ and $\mu_T$ witness the left fair amenability of $S$ and $T$ respectively. Let $\pi_S, \pi_T$ denote the projections from $\Sc{P}(S\times T)$ onto $\Sc{P}(S)$ and $\Sc{P}(T)$, respectively.
\begin{pruf}
\step{define}{Define $\mu$, on $S\times T$, for each rectangle $R = A\times B$ where $A\subseteq S$ and $B\subseteq T$:
\[
	\fn\mu{R} \defeq \fn{\mu_S}{\fn{\pi_S}R}\fn{\mu_T}{\fn{\pi_T}R} = \fn{\mu_S}{A}\fn{\mu_T}{B},
\] which, while not yet defined for all subsets of $S\times T$, is clearly left fairly invariant and finitely-additive, and with $\fn\mu{S\times T} = \fn{\mu_S}{S}\fn{\mu_T}{T} = 1$.}
\step{extend1}{It follows that
\[
    \fn\mu{\bigcup_{i\in I} R_i} = \sum_{i\in I} \fn\mu{R_i},
\] for each finite collection of disjoint rectangles\footnote{Take care to avoid confusing \emph{finite collections of rectangles} with \emph{collections of finite rectangles}.} $\setof{R_i}_{i\in I}$, and this is also left fairly invariant.}
    \begin{pruf}
    \pf\ If $(s,t)$ acts injectively on $\bigcup_{i\in I} R_i$, then $s$ acts on $\fn{\pi_S}{R_i}$ injectively for each $i\in I$, likewise for $t\in T$ on $\fn{\pi_T}{R_i}$. Furthermore, $(s,t)$ preserves the disjointness of $\setof{R_i}_{i\in I}$. 
    \end{pruf}
\step{extend2}{Let $C$ be an arbitrary subset of $S\times T$. $C$ is not necessarily a rectangle, so extend $\mu$ using 
\[
	\fn\mu{C} \defeq \sup \fn\mu{\bigcup_{i\in I} R_i},
\] where the supremum is taken over all finite collections of subrectangles of $C$.}
\step{verify}{$\mu$ is then defined for all subsets $C$ of $S\times T$, and is left fairly invariant.}
    \begin{pruf}
    \pf\ If $(s,t)\in S\times T$ acts injectively on $C$ then it acts injectively on any finite collection of disjoint subrectangles of $C$. Each finite collection of disjoint subrectangles of $(s,t)C$ has the form $\setof{(s,t)R_i}_{i\in I}$ for a finite collection of disjoint subrectangles $\setof{R_i}_{i\in I}$ of $C$. Hence 
    \begin{align*}
        \fn\mu{(s,t)C} &= \sup\fn\mu{\bigcup_{i\in I} (s,t)R_i} \\
            &= \sup\sum_{i\in I} \fn\mu{(s,t) R_i} \\
            &= \sup\sum_{i\in I} \fn\mu{R_i} \\
            &= \sup\fn\mu{\bigcup_{i\in I} R_i} \\
            &= \fn\mu{C},
    \end{align*}
     as required. \qed
    \end{pruf}
\end{pruf}
}
\end{theorem}

\anondivision
Another result for groups states that every left amenable group is also right amenable, and furthermore, a left invariant measure and right invariant measure can be combined to provide a bi-invariant measure \citep[p148]{Wagon:1993th}. This result doesn't hold for all semigroups (either classically or fairly), but a similar proof technique can be applied to the fair amenability of semigroups with involution.

\begin{lemma}\label{lem:involutionleftright}
Let $S$ be a semigroup with involution $\ast$. If $S$ is left fairly amenable, then it is right fairly amenable (and vice-versa). 
\Proof{
$A^\ast := \setof{a^\ast: a\in A}$, and so $(As)^\ast = s^\ast A^\ast$. Suppose that $\mu$ on $S$ is left fairly invariant and define $\nu$ by setting $\fn\nu{A} = \fn\mu{A^\ast}$ for all $A$. 

\begin{pruf}
\step{fapm}{$\nu$ has total measure 1.}
    \begin{pruf}
    \pf\ 
    \begin{align*}
        \fn\nu{S} &= \fn\nu{S^\ast} \\
            &= \fn\mu{S} = 1.
    \end{align*}
    \end{pruf}
\step{fa}{$\nu$ is finitely additive.}
    \begin{pruf}
    \pf\ For all disjoint $A,B \subseteq S$, 
    \begin{align*}
        \fn\nu{A\cup B} &= \fn\mu{(A\cup B)^\ast} \\
            &= \fn\mu{A^\ast \cup B^\ast} \\
            &= \fn\mu{A^\ast} + \fn\mu{B^\ast} \\
            &= \fn\nu{A} + \fn\nu{B}.
    \end{align*}
    \end{pruf}
\step{rightinvar}{$\nu$ is right fairly invariant.}
    \begin{pruf}
    \pf\ If $s$ acts injectively on the \emph{right} of $A$, then for $a,b\in A$,
\begin{align*}
	s^\ast a^\ast = s^\ast b^\ast &\iff (as)^\ast = (bs)^\ast\\
		&\iff as = bs \\
		&\imp a = b \\
		&\iff a^\ast = b^\ast
\end{align*}
and so $s^\ast$ acts injectively on the \emph{left} of $A^\ast$. Then
\[
	\fn\nu{As} = \fn\mu{s^\ast A^\ast} = \fn\mu{A^\ast} = \fn\nu{A}
\] wherever $s$ acts injectively on the right of $A$. \qed
    \end{pruf}
\end{pruf}
}
\end{lemma}

Thus groups, inverse semigroups, semigroups of binary relations, and all other $\ast$-semigroups join the commutative semigroups as classes of semigroups where each example is either \emph{fairly amenable (both ways)}, or not at all.

In the next section I give an example of a semigroup that is fairly amenable on one side but not the other.

Every subgroup of an amenable group is amenable, including those subgroups having measure zero. A quick summary of this proof goes as follows: let $G$ be an amenable group with measure $\mu$, $H$ a subgroup. Choose a set $M$ of representatives from each left coset of $H$, then define a measure $\nu$ on $H$ by setting $\fn\nu{A} \defeq \fn\mu{MA}$ for all $A\subseteq H$ \citep[p149]{Wagon:1993th}. It would be nice to emulate this in the semigroup case, but it seems there is no adequate analogue for semigroups of the coset structure of a group. Perhaps the obvious should be stated:

\begin{lemma}\label{lem:subsgpinheritFAifmass}
Let $S$ be a left [right] fairly amenable semigroup with measure $\mu$, and let $T$ be a subsemigroup of $S$ having $\fn\mu{T} > 0$. $T$ is then left [right] fairly amenable.
\Proof{
We may use $\nu$ as given by $\fn\nu{A} = \fn\mu{A} / \fn\mu{T}$ for all $A\subseteq T$.
\qed}
\end{lemma}

This mirrors the classical case \citep[p.518]{Day1956}. In particular, any subgroup $G$ of a left or right fairly amenable semigroup is amenable \emph{provided that $\fn\mu{G} > 0$}. 

\begin{corollary}\label{cor:addingazero}
Let $S$ be a semigroup without zero. $S^0$ is left [right] fairly amenable if and only if $S$ is.
In particular, if $G$ is a group, $G^0$ is fairly amenable if and only if $G$ is amenable.
\Proof{
Since the finite case is trivial , assume that $S$ is infinite. If $S^{0}$ is left fairly amenable with $\mu'$, since $S^{0}$ contains a zero, by Lemma \ref{lem:finitemassless} $\fn{\mu'}{\setof{0}} = 0$, which by finite additivity implies $\fn{\mu'}{S} = 1$. By Lemma \ref{lem:subsgpinheritFAifmass} $S$ is fairly amenable and, in the case of a group, amenable by Corollary \ref{cor:backwardscompatible}. 

Conversely, if $S$ is left fairly amenable with some $\mu$ then assigning $\fn{\mu'}{A} = \fn\mu{A\cap S}$ yields a fairly invariant measure $\mu'$ on $S^0$. The case on the right holds similarly.\qed
}
\end{corollary}

0-groups are examples of Clifford semigroups, which in turn are charactarised as being strong semilattices of groups \citep[p94]{Howie1976}, and in turn are examples of inverse semigroups. One wonders, therefore, what we can say about Clifford semigroups in general. The following example furnishes us with both a fairly amenable Clifford semigroup that is not a 0-group, having a non-amenable subgroup in a non-trivial manner.

\begin{example}\label{ex:cliffordF2andF1}
Let $S$ be the union of two groups as follows: set $G \cong \mathbb{F}_2$ (not amenable) and $H \cong \mathbb{F}_1$ (amenable), and let $\phi: G\to H$ be the homomorphism mapping $x\mapsto 1_H$ for all $x\in G$. Define the operation on $S$ as a strong semilattice $Y = (\setof{1,0},\wedge)$ of the groups $G,H$, i.e. if one of $x$ or $y$ is in $H$ we map the other via $\phi$ into $H$ to compute $xy$. Despite the presence of $\mathbb{F}_2$, $S$ \emph{is} fairly amenable.
\Proof{
Let $\mu_H$ witness the amenability of $H$. Define for $S$ the measure $\mu$ given by
\[
	\fn\mu{A} \defeq \fn{\mu_H}{\fn\phi{A\cap G} \cup (A\cap H)} \quad\text{for all } A \subseteq S,
\]
which is invariant under action of $H$. Since $H$ is an infinite $\Sc{H}$-class, $\fn{\mu_H}{\setof{1_H}} = 0$ by Lemma \ref{lem:finitemassless2}, and therefore $\fn\mu{G} = 0$; it follows that $\fn\mu{A} = \fn{\mu_H}{A\cap H}$ for any $A\subseteq S$. If $A\subseteq H$ then $gA = A = Ag$ for all $g\in G$, so $\mu$ is trivially invariant under $G$, and thus $\mu$ suffices. \qed
}
\end{example}

The following example shows a fairly amenable Clifford semigroup that has no amenable subgroup as part of the semilattice.

\begin{example} Consider the semilattice on the integers $Y = (\mathbb{Z},\wedge)$ where $a\wedge b = \min\setof{a,b}$ for all $a,b\in Y$, together with a measure $\mu$ derived from the \Folner sequence given by $F_n = [-n,n] \cap Y$ for each $n$. 

Now $\fn\mu{k\wedge Y} = \fn\mu{(-\infty,k] \cap Y} = \frac12$ for all $k\in Y$, all finite sets have measure $0$, and the semilattice is fairly amenable. 

Suppose we take $S$ to be a strong semilattice of infinitely many non-amenable groups, as follows:
\begin{itemize}
\item Let the semilattice $Y$ be isomorphic to $(\mathbb{Z},\wedge)$, as previously;
\item For each $k\in\mathbb{Z}$ let $G_k$ be a non-amenable group;
\item For each $k\in\mathbb{Z}$ let $\nu_k$ be any finitely-additive measure on $G_k$ with $\fn{\nu_k}{G_k} = 1$ (which is necessarily not invariant).
\end{itemize}

We can extend the $\mu$ given on $Y$ to a fairly-invariant $\mu_S$ on $S$ by setting, for a fixed free ultrafilter $U$ over $\mathbb{N}$,
\[
	\fn{\mu_S}{A} = \lim_U \frac1{2n+1} \sum_{k=-n}^n \fn{\nu_k}{G_k\cap A}.
\]
While every $G_k$ is not amenable, $\mu_S$ witnesses the fair amenability of $S$. \qed
\end{example}

\begin{corollary}
If the Clifford semigroup $S$ is a strong \emph{finite} semilattice $Y$ of groups and $S$ is fairly amenable, at least one of the groups is amenable.
\Proof{
Suppose all the groups in $\setof{G_y: y\in Y}$ are non-amenable, and the finitely-additive measure $\mu$ witnesses the fair amenability of $S$.
\begin{pruf}
\step{eachmassiszero}{$\fn\mu{G_y} = 0$ for all $y\in Y$.}
    \begin{pruf}
    \pf\ Use Lemma \ref{lem:subsgpinheritFAifmass}. 
    \end{pruf}
\step{allmassiszero}{$1 = \fn\mu{S} = 0$, contradiction.}
    \begin{pruf}
    \pf\ $S=\bigcup_{y\in Y}G_y$, which is a disjoint union, and then as there are only finitely many groups in the semilattice, $\fn\mu{S} = 0$. \qed
    \end{pruf}
\end{pruf}
}
\end{corollary}

\anondivision
One final theorem on groups that translates well to fairly amenable semigroups is that a directed union of amenable groups is also amenable. 
\begin{theorem}\label{thm:directunionoffairlyamenablesemigroups}
If $S$ is the directed union of left [right] fairly amenable semigroups, then $S$ is left [right] fairly amenable.
\Proof{
This proof uses essentially the same topological argument as in \cite[p150]{Wagon:1993th}.
Let $\setof{S_i: i\in I}$ be the directed system of left fairly amenable semigroups whose union is $S$: i.e. for each $a,b\in I$ there exists a $c\in I$ such that $S_a$ and $S_b$ are subsemigroups of $S_c$, and, $S = \bigcup_{i\in I} S_i$. For each $i\in I$:
\begin{itemize} 
\item let $\mu_i$ be the left fairly invariant finitely-additive measure corresponding to $S_i$;
\item let $M_i$ be the set of finitely-additive measures $m:\Sc{P}(S)\to[0,1]$ such that $\fn{m}{S} = 1$ and whenever $s\in S_i$ acts injectively on $A\subseteq S$, $\fn{m}{sA} = \fn{m}{A}$.
\end{itemize}

\begin{pruf}
\step{nonempty}{$M_i$ is non-empty for all $i\in I$.}
    \begin{pruf}
    \pf\ Define $\fn{m_i}{A} \defeq \fn{\mu_i}{A\cap S_i}$ for all $A\subseteq S$. Clearly $m_i\in M_i$.
    \end{pruf}
\step{closed}{Each $M_i$ is a closed subset of $[0,1]^{\Sc{P}(S)}$.}
    \begin{pruf}
    \pf\ Suppose $f\notin M_i$; either $f$ fails to be finitely additive, fails to be left fairly invariant for some $s\in S_i$, or $f(S) \ne 1$. It is possible to vary the ``amount'' by which each of the three conditions is violated (e.g. $1 - f(S) = \epsilon$), thus forming an open neighborhood of $f$ consisting of points behaving similarly. This argument is similar to \cite[p126]{Wagon:1993th}.
    \end{pruf}
\step{finiteintersection}{The collection $\setof{M_i: i\in I}$ has the finite intersection property.}
    \begin{pruf}
    \pf\ If $S_a, S_b \subseteq S_c$ then $M_a\cap M_b \supseteq M_c$, since each member must be left fairly invariant for increasingly many elements. 
    \end{pruf}
\step{compactness}{There exists some $\mu\in\bigcap_{i\in I} M_i$ which is the required left fairly-invariant measure.}
    \begin{pruf}
    \pf\ From Tychonoff's Theorem, the space $[0,1]^{\Sc{P}(S)}$ is compact; equivalently, any collection of closed subsets with the finite intersection property is nonempty, and $\setof{M_i: i\in I}$ is an example of such a collection. 
    \end{pruf}
The right case is handled analogously. \qed
\end{pruf}
}
\end{theorem}


\section{Examples}

\begin{proposition}\label{prop:freecommutativesgp}
Any finitely-generated free Abelian semigroup, such as $(\mathbb{N},+)$, is fairly amenable.
\Proof{
The free Abelian semigroup on $k$ generators is isomorphic to $(\mathbb{N}\cup\setof{0})^k$ minus the origin, and again every action is injective. The \Folner sequence given by \(
    F_n = \setof{(a_1, a_2, \ldots, a_k): a_1, a_2, \ldots a_k < n}
\) suffices.
\qed
}
\end{proposition}

\begin{proposition}
$(\mathbb{N},\cdot)$ (the natural numbers with multiplication) is also a cancellative Abelian semigroup. However, it is infinitely generated (by the primes). It is also fairly amenable. 
\Proof{As usual a totally invariant finitely-additive measure is required. There exists a \Folner sequence $\setof{F_n}_{n\in\mathbb{N}}$ where $F_n$ consists of the products of powers of the first $n$ primes, and each power lies in $[0,n]$, i.e.  
\[
    F_n \defeq \setof{p_1^{i_1} p_2^{i_2}\cdots p_n^{i_n}: 0\le i_j \le n, j=1, \ldots, n},
\]
as required. \cite{Bergelson:2005tx} demonstrated a family of \Folner sequences of this kind. \qed 
}
\end{proposition}

\begin{example}\label{ex:freesemigroup2gens}
The free semigroup on two generators $FS_2 = \setof{a,b}^{+}$ is neither left nor right fairly amenable.
\Proof{
Suppose $S = \setof{a,b}^{+}$ is left fairly amenable and $\mu$ is the required measure. Note that $a$ and $b$ both act injectively on $S$\ and so we require $\fn\mu{aS} = \fn\mu{S} = \fn\mu{bS}$. 
But since $S = \setof{a,b} \dotcup aS \dotcup bS$, 
\[
	1 = \fn\mu{S} = \fn\mu{\setof{a,b}} + \fn\mu{aS} + \fn\mu{bS} = \fn\mu{\setof{a,b}} + 1 + 1 \ge 2,
\]
contradiction. By a similar argument, $FS_2$ is not right fairly amenable. (Alternatively, endow the semigroup with an involution $*$ where $a^* \defeq b$ and vice-versa, and apply Lemma \ref{lem:involutionleftright}.)
\qed}

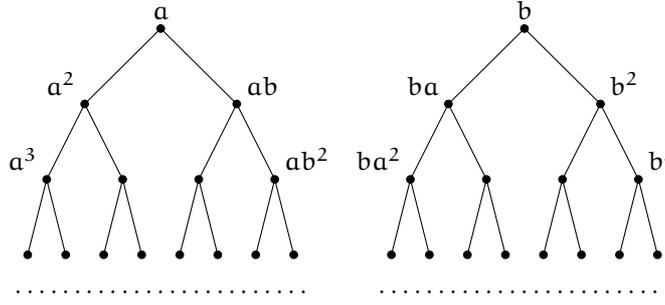
\begin{figure}
\begin{center}
\begin{tikzpicture}
\filldraw[fill=black] (0,3) circle (0.05) node[anchor=south] {$a$}
                      (-1,2) circle (0.05) node[anchor=south east] {$a^2$}
                      (1,2) circle (0.05) node[anchor=south west] {$ab$}
                      (-1.5,1) circle (0.05) node[anchor=south east] {$a^3$}
                      (-0.5,1) circle (0.05) 
                      (0.5,1) circle (0.05)
                      (1.5,1) circle (0.05) node[anchor=south west] {$ab^2$}
                      (-1.75,0) circle (0.05)
                      (-1.25,0) circle (0.05)
                      (-0.75,0) circle (0.05)
                      (-0.25,0) circle (0.05)
                      (0.25,0) circle (0.05)
                      (0.75,0) circle (0.05)
                      (1.25,0) circle (0.05)
                      (1.75,0) circle (0.05);
\draw (-1,2) -- (0,3) -- (1,2);
\draw (-1.5,1) -- (-1,2) -- (-0.5,1);
\draw (1.5,1) -- (1,2) -- (0.5,1);
\draw (-1.75,0) -- (-1.5,1) -- (-1.25,0);
\draw (-0.75,0) -- (-0.5,1) -- (-0.25,0);
\draw (1.75,0) -- (1.5,1) -- (1.25,0);
\draw (0.75,0) -- (0.5,1) -- (0.25,0);
\draw (0,-0.5) node {$\cdots\cdots\cdots\cdots\cdots\cdots\cdots\cdots\cdot$};
\end{tikzpicture}
\begin{tikzpicture}
\filldraw[fill=black] (0,3) circle (0.05) node[anchor=south] {$b$}
                      (-1,2) circle (0.05) node[anchor=south east] {$ba$}
                      (1,2) circle (0.05) node[anchor=south west] {$b^2$}
                      (-1.5,1) circle (0.05) node[anchor=south east] {$ba^2$}
                      (-0.5,1) circle (0.05)
                      (0.5,1) circle (0.05) 
                      (1.5,1) circle (0.05) node[anchor=south west] {$b^3$}
                      (-1.75,0) circle (0.05)
                      (-1.25,0) circle (0.05)
                      (-0.75,0) circle (0.05)
                      (-0.25,0) circle (0.05)
                      (0.25,0) circle (0.05)
                      (0.75,0) circle (0.05)
                      (1.25,0) circle (0.05)
                      (1.75,0) circle (0.05);
\draw (-1,2) -- (0,3) -- (1,2);
\draw (-1.5,1) -- (-1,2) -- (-0.5,1);
\draw (1.5,1) -- (1,2) -- (0.5,1);
\draw (-1.75,0) -- (-1.5,1) -- (-1.25,0);
\draw (-0.75,0) -- (-0.5,1) -- (-0.25,0);
\draw (1.75,0) -- (1.5,1) -- (1.25,0);
\draw (0.75,0) -- (0.5,1) -- (0.25,0);
\draw (0,-0.5) node {$\cdots\cdots\cdots\cdots\cdots\cdots\cdots\cdots\cdot$};
\end{tikzpicture}
\caption{The right Cayley graph for the free semigroup on two generators $\setof{a,b}^+$.}
\end{center}
\end{figure}
\end{example}

\begin{remark}
Note that the previous argument can be adapted to any finite number of generators $n\ge 2$. Note also that $FS_2^0$ (the free semigroup on two generators with a zero adjoined) is now not fairly amenable either, in contrast to the classical case.
\end{remark}

\begin{remark}
Another theorem on groups states that if a group $G$ is amenable and $N\lhd G$, then $G/N$ is also amenable; since every congruence on a group arises as the cosets of a normal subgroup this means that every quotient of an ameable group is amenable. Given $\mu$ on an amenable $G$ we may set $\nu$ on $G/N$ using
\begin{equation}\label{eqn:quotientmeasure}
	\fn\nu{A} = \fn\mu{\bigcup A}.
\end{equation}

The corresponding situation in fairly amenable semigroups encounters problems. Let $\sigma$ be a congruence on a fairly left amenable semigroup $S$ with measure $\mu$. Clearly $\nu$ has total measure 1 and is finitely-additive. However it is not always going to be left fairly invariant. 
\end{remark}

\begin{example}
As described in Proposition \ref{prop:freecommutativesgp}, the free Abelian semigroup on two generators $S$ is fairly amenable with the measure $\mu$. Let $\sigma$ be the congruence on $S$ with $(b,b^2), (b, ab)\in \sigma$, i.e.
\[
	S/\sigma \cong \mathrm{sgp}\pres{a,b}{ab=ba=b^2=b}.
\]
$S/\sigma$ is fairly amenable (it is a free commutative semigroup on one generator with a zero), however $\nu$ as in Equation \ref{eqn:quotientmeasure} is not fairly invariant since
\(
	\fn\nu{A} = \fn\nu{(b\sigma)^{-1}A}
\) (the Dirac delta measure), via Lemma \ref{lem:finitemassless2}.	 \qed
\end{example}

Now we consider some bands. Recall that, in the classical theory, a right zero semigroup is left amenable but not right amenable. 

\begin{example}\label{ex:zerosemigroups}
Let $S$ be a left (or right) zero semigroup. $S$ is fairly amenable (both sides).
\Proof{
The finite case is handled by Corollary \ref{cor:allfinite}, so assume $S$ is an infinite left zero semigroup. 

\begin{pruf}
\step{right}{Any finitely-additive measure $\mu$ with $\mu(S) = 1$ is right fairly invariant.}
    \begin{pruf}
    \pf\ For any $A\subseteq S$ and $s\in S$, $As = A$, so $\fn\mu{As} = \fn\mu{A}$ trivially.
    \end{pruf}
\step{left}{There are infinitely many finitely-additive measures $\mu$ with $\fn\mu{S} =1$ that are left fairly invariant.} 
    \begin{pruf}
    \pf\ For any $A\subseteq S$ and $s\in S$, $sA = \setof{s}$, and by Lemma \ref{lem:finitemassless} every $\fn\mu{\setof{s}} = 0$ if $\mu$ is fairly invariant, but since singletons are the only sets injectively acted on on the left, the following suffices. Fix any free ultrafilter $U$, and define $\fn\mu{A} = \fn{\chi_U}{A}$. 
    \end{pruf}
The argument holds on the right analogously. \qed
\end{pruf}}
\end{example}

\begin{example}
Every rectangular band is fairly amenable.
\Proof{We have just seen the specific examples of left and right zero semigroups (Example \ref{ex:zerosemigroups}). Each rectangular band is isomorphic to the product of a left zero semigroup and a right zero semigroup, therefore by Theorem \ref{thm:directproductfairly} all rectangular bands are fairly amenable.\qed}
\end{example}

\begin{figure}
\begin{center}
\begin{tikzpicture}
\draw[gray!20,->] (0,0) -- (0,2);
\draw[gray!20,->] (0,0) -- (2,0);
\filldraw[fill=gray!20,draw=black,very thick] (0,0) rectangle (1,1);
\draw (0.5,-0.05) node[anchor=north] {$\square_n$};
\end{tikzpicture}
\begin{tikzpicture}
\draw[gray!20,->] (0,0) -- (0,2);
\draw[gray!20,->] (0,0) -- (2,0);
\draw[gray!50] (0,0) rectangle (1,1);
\filldraw[fill=gray!20,draw=black,very thick] (0,0) rectangle (0.7,1);
\draw[very thick] (0,0) -- (0, 1.3);
\draw (0.35,-0.05) node[anchor=north] {$p^k\square_n$};
\end{tikzpicture}
\begin{tikzpicture}
\draw[gray!20,->] (0,0) -- (0,2);
\draw[gray!20,->] (0,0) -- (2,0);
\draw[gray!50] (0,0) rectangle (1,1);
\filldraw[fill=gray!20,draw=black,very thick] (0.5,0) rectangle (1.2,1);
\draw[very thick] (0.5,0) -- (0.5, 1.3);
\draw (0.9,-0.05) node[anchor=north] {$q^jp^k\square_n$};
\end{tikzpicture}
\begin{tikzpicture}
\draw[gray!20,->] (0,0) -- (0,2);
\draw[gray!20,->] (0,0) -- (2,0);
\draw[gray!50] (0,0) rectangle (1,1);
\filldraw[fill=gray!20,draw=black,very thick] (0,0) rectangle (0.5,1);
\filldraw[fill=gray!20,draw=black,very thick] (1,0) rectangle (1.2,1);
\draw[very thick] (0.5,1) -- (0.5, 1.3);
\draw (1,-0.05) node[anchor=north] {$q^jp^k\square_n\,\triangle\,\square_n$};
\end{tikzpicture}
\caption{\label{fig:sizeofppsquare} Deriving $\abs{q^jp^k\square_n\,\triangle\,\square_n}$ in the bicyclic monoid.}
\end{center}

\end{figure}
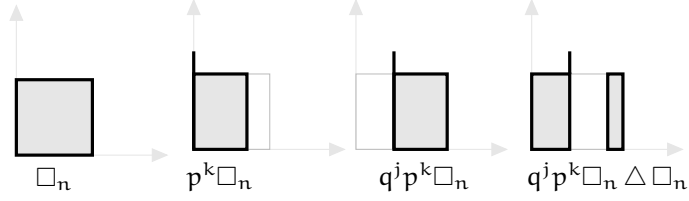

\begin{example}\label{ex:bicyclicmonoid}
The bicyclic monoid $B$ is fairly amenable.
\Proof{
Recall that $B = \text{mon}\pres{p,q}{pq = 1} = \setof{q^mp^n: m,n \in \mathbb{N} \cup\setof{0}}$. 
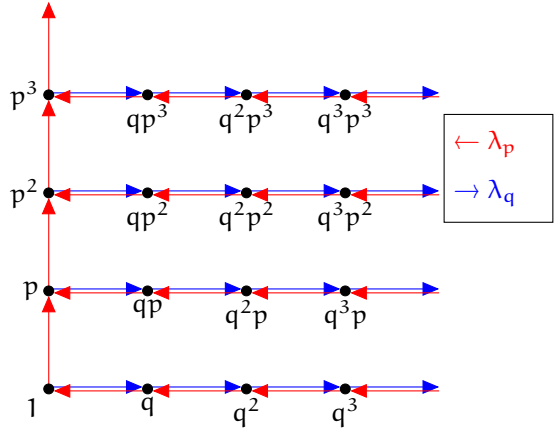
\begin{figure}
\begin{center}
\begin{tikzpicture}[scale=1.3]
\foreach \y in {0,...,3}
{
    \pgfmathadd{\y}{0.95}\pgfmathsetmacro{\yp}{\pgfmathresult};
    \draw[red,->] (0,\y) -- (0,\yp);
}
\foreach \x in {0,...,3}
{
    \pgfmathadd{\x}{0.95}\pgfmathsetmacro{\xp}{\pgfmathresult};
    \pgfmathadd{\x}{0.05}\pgfmathsetmacro{\xq}{\pgfmathresult};
    \foreach \y in {0,...,3}
    {
        \pgfmathadd{\y}{0.02}\pgfmathsetmacro{\yp}{\pgfmathresult};
        \pgfmathadd{\y}{-0.02}\pgfmathsetmacro{\yq}{\pgfmathresult};
        \draw[blue,->] (\x,\yp) -- (\xp,\yp);
        \draw[red,<-] (\xq,\yq) -- (\xp,\yq);
        \filldraw[fill=black] (\x,\y) circle (0.05);
    }
}
\draw (0,0) node[anchor=north east]{$1$};
\draw (0,1) node[anchor=east]{$p$};
\draw (0,2) node[anchor=east]{$p^2$};
\draw (0,3) node[anchor=east]{$p^3$};
\draw (1,0) node[anchor=north] {$q$};
\draw (1,1) node[anchor=north] {$qp$};
\draw (1,2) node[anchor=north] {$qp^2$};
\draw (1,3) node[anchor=north] {$qp^3$};
\foreach \x in {2,3}
{
    \draw (\x,0) node[anchor=north] {$q^{\x}$};
    \draw (\x,1) node[anchor=north] {$q^{\x}p$};
    \draw (\x,2) node[anchor=north] {$q^{\x}p^2$};
    \draw (\x,3) node[anchor=north] {$q^{\x}p^3$};
}
\draw (4,1.7) rectangle (5.1,2.8);
\draw (4,2.5) node[anchor=west,red]{$\leftarrow \lambda_p$};
\draw (4,2) node[anchor=west,blue]{$\rightarrow \lambda_q$};
\end{tikzpicture}
\caption{\label{fig:bicyclicgraph}Part of the left Cayley graph of the bicyclic monoid $B$.}
\end{center}
\end{figure}

Consider the sequence given by $\square_n = \setof{q^jp^k: j, k \le n}$ for all $n\in \mathbb{N}$. It will suffice to show this sequence is \Folner for any element on the left. 

The element $q$ acts injectively on the left of all $B$, so $\abs{q\square_n} = \abs{\square_n}$ and $\abs{q\square_n\,\triangle\,\square_n} = 2n$. $p$ on the other hand does not act injectively on the left of $\square_n$, in which case $\abs{p\square_n} \le \abs{\square_n}$. Since the minimal non-injective sets for each left multiplication by $p$ are $\setof{p^k, qp^{k+1}}$ for each $k$, we can see exactly that $\abs{p\square_n} = (n-1)n + 1$, and $\abs{p\square_n\,\triangle\,\square_n} = n+1$. For any arbitrary $x = q^jp^k$, then,
\[
    \abs{x\square_n\,\triangle\,\square_n} = k + n(2j-k) \quad\text{for all } n > j
\]
(depicted in Figure \ref{fig:sizeofppsquare}) which is linear in $n$, and therefore the \Folner sequence $\setof{\square_n}_{n\in\mathbb{N}}$ suffices.

$B$ is inverse, so Lemma \ref{lem:involutionleftright} applies and hence $B$ is fairly amenable on both sides. \qed}
\end{example}

\begin{example}\label{ex:polycyclicmonoid}
The polycyclic monoid on two generators, $P_2$, is not fairly amenable. As described by \cite{Milan:2008tp}, $P_2$ has the weak containment property, so it follows that fair amenability is not equivalent to weak containment.
\Proof{
Recall that 
\[
	P_2 = \text{mon}^0\pres{p,q,p^{-1},q^{-1}}{pp^{-1}=1=qq^{-1}, pq^{-1}=0=qp^{-1}},
\] 
and so every element other than 0 or 1 can be written canonically in the form $x^{-1}y$, where $x,y$ are (possibly empty) strings over the alphabet $\setof{p,q}$ \citep{Lawson2004}. It follows that (at least) the free monoids $\setof{p^{-1},q^{-1}}^{*}$ and $ \setof{p,q}^{*}$ are embedded within $P_2$. 

\begin{pruf}
\step{setup}{Assume $P_2$ is left fairly amenable with measure $\mu$, and for each $x\in P_2$ let $H_x\subseteq P_2$ consist of elements with their canonical form starting with the string $x$. $P_2$ can be decomposed like so: 
\[
	P_2 = H_{p^{-1}} \dotcup H_{q^{-1}} \dotcup H_p \dotcup H_q \dotcup \setof{0, 1}. 
\]}
\step{injectives}{Consider the injective left actions $\lambda_{p^{-1}}, \lambda_{q^{-1}}$
\[
	p^{-1} P_2 = H_{p^{-1}} \dotcup \setof{0}, \quad q^{-1} P_2 = H_{q^{-1}} \dotcup \setof{0}.
\]}
\step{apply}{Apply $\mu$ to see that it is not left fairly invariant.}
    \begin{pruf}
    \pf\ \begin{align*}
        	1 &= \fn\mu{P_2} \\
        	&= \fn\mu{ H_{p^{-1}} \dotcup H_{q^{-1}} \dotcup H_p \dotcup H_q \dotcup \setof{0, 1}} \quad\because\text{step \stepref{setup}} \\
        	&= \fn\mu{H_{p^{-1}}} + \fn\mu{H_{q^{-1}}} + \fn\mu{H_p} + \fn\mu{H_q} + \fn\mu{\setof{0, 1}} \\
        	&= \fn\mu{H_{p^{-1}}} + \fn\mu{H_{q^{-1}}} + \fn\mu{H_p} + \fn\mu{H_q} \quad\because\text{Lemma \ref{lem:finitemassless}}\\
        	&= \fn\mu{p^{-1}P_2} + \fn\mu{q^{-1}P_2} + \fn\mu{H_p} + \fn\mu{H_q} \quad\because\text{step \stepref{injectives}} \\
        	&= 1 + 1 + \fn\mu{H_p} + \fn\mu{H_q} \quad\because\text{fair invariance} \\
        	&\ge 2,
        \end{align*}
        contradiction. 
    \end{pruf}
\end{pruf}
$P_2$ is also inverse, so by Lemma \ref{lem:involutionleftright} it is not right fairly amenable either. 
\qed}
\end{example}

\begin{remark}As with $FS_2$ and greater, the previous argument can be adapted to any finite number of generators $n\ge 2$. $P_2$ is also an example of an inverse semigroup that is not fairly amenable, but is classically amenable because the maximal group homomorphic image (the trivial group) is amenable.
\end{remark}

\begin{example}\label{ex:baerlevi}
For a Levi-Baer semigroup $LB(p,q)$,
\begin{enumerate}[(i)]
\item $LB(p,q)$ is not left fairly amenable if $p=q$;
\item $LB(p,q)$ is not right fairly amenable for all $p,q$.
\end{enumerate}
\Proof{
Recall that Levi-Baer semigroups are left cancellative, left simple, and have no idempotents.\footnote{A Baer-Levi semigroup $BL(p,q)$ is defined as being the set of injective maps $f$ on some infinte set $X$ having cardinality $p$, such that $\abs{X\backslash \fn{f}{X}}$ is some fixed infinite cardinal $q\le \abs{X}=p$ \citep{clifford1967algebraic}. Conventionally, products in Baer-Levi semigroups are written in ``algebraist'' order---the composition of $f$ and $g$ is written $fg$---and hence the Baer-Levi semigroups are normally referred to as right cancellative and right simple. However, to remain consistent, I shall use $\circ$ and consider the equivalent ``Levi-Baer'' semigroup, which is left cancellative and left simple.} For succinctness let $S$ be shorthand for $LB(p,q)$.

On the left:
\begin{pruf}
\step{setup}{Let $a,b\in S$ be such that the right ideals $a\circ S$ and $b\circ S$ are disjoint. (There are two disjoint right ideals if, and only if, $p=q$.) For example, if $S$ is the Baer-Levi semigroup on $\mathbb{N}$, we may pick $a:n\mapsto 2n$ and $b:n\mapsto 2n+1$. Let $R = S\backslash ((a\circ S) \cup (b\circ S))$.}
\step{injective}{Since $S$ is left cancellative, every left action is injective.}
\step{apply}{Assume $S$ is left fairly amenable with measure $\mu$, and derive a contradiction.}
    \begin{pruf}
    \pf\ \begin{align*}
    	1 &= \fn\mu{S} \\
    	  &= \fn\mu{(a\circ S) \dotcup (b\circ S) \dotcup R} \quad\text{by defintion}\\
    	  &= \fn\mu{a\circ S} + \fn\mu{b\circ S} + \fn\mu{R} \\
    	  &= \fn\mu{S} + \fn\mu{S} + \fn\mu{R} \quad\because\text{left fairly invariant}\\
    	  &\ge 2,
    \end{align*}
    a clear contradiction. 
    \end{pruf}
\end{pruf}

On the right: 
\begin{pruf}
\step{setup}{For each $s\in S$ let the equivalence relation $\theta_s$ be given by $a\,\theta_s\,b \iff a\circ s = b\circ s$ for all $a,b\in S$. Since $S$ consists of maps on the set $X$, $\theta_s$ depends only on $\fn{s}{X}$, so $a\,\theta_s\,b \iff a|_{\fn{s}{X}} = b|_{\fn{s}{X}}$.}
    \begin{pruf}
    \pf\ For any $a,b,s\in S$,
        \begin{align*}
            a\,\theta_s\,b &\iff a\circ s = b\circ s \\
                &\iff \fn{a}{\fn{s}{x}} = \fn{b}{\fn{s}{x}} \text{ for all } x\in X \\
                &\iff \fn{a}{y} = \fn{b}{y} \text{ for all } y\in \fn{s}{X} \\
                &\iff a|_{\fn{s}{X}} = b|_{\fn{s}{X}}.
        \end{align*}
    \end{pruf}
\step{infinite}{For every $s\in S$, every $\theta_s$-equivalence class is nonempty and infinite.}
    \begin{pruf}
    \pf\ By definition $\abs{X\backslash\fn{s}{X}}$ is some infinite cardinal, therefore a Baer-Levi semigroup on $X\backslash\fn{s}{X}$ can be used to generate elements of each $\theta_s$-class.
    \end{pruf}
\step{what}{For each $s\in S$ define two disjoint sets $M_1, M_2$ by choosing two distinct elements from each $\theta_s$-class. $S\circ s = M_1\circ s = M_2\circ s$ and while the action $S\circ s$ is not injective, the actions on $M_1$ and $M_2$ are injective.}
    \begin{pruf}
    \pf\ By definition, $\theta_s$ partitions $S$ into sets that map to the same element under the right action of $s$, so $S\circ s = M_1\circ s = M_2\circ s$. For any $a,b\in M_1$, $a\circ s = b\circ s\imp a\,\theta_s\,b \imp a=b$, similarly for $M_2$. 
    \end{pruf}
\step{contradict}{Assume that $S$ is right fairly amenable with measure $\nu$. This results in a contradiction.}
    \begin{pruf}
    \pf\ Let $R = S\backslash (M_1\cup M_2)$, then
    \begin{align*}
        1 &= \fn\nu{S} \\
            &= \fn\nu{M_1 \dotcup M_2 \dotcup R} \quad\because\text{definition} \\
            &= \fn\nu{M_1} + \fn\nu{M_2} + \fn\nu{R} \\
            &= \fn\nu{M_1\circ s} + \fn\nu{M_2\circ s} + \fn\nu{R} \quad\because\text{right fairly invariant} \\
            &= \fn\nu{S\circ s} + \fn\nu{S\circ s} + \fn\nu{R} \quad\because\text{step \stepref{what}} \\
            &= \fn\nu{S} + \fn\nu{S} + \fn\nu{R}  \quad\because S\text{ is left simple} \\
            &= 1 + 1 + \fn\nu{R} \\
            &\ge 2, 
    \end{align*}
    a clear contradiction. \qed
    \end{pruf}
\end{pruf}
} 
\end{example}

\begin{example}\label{ex:leftbutnotright}
Left groups are left simple, right cancellative semigroups that are characterised as being direct products of groups and left zero semigroups. Let $Z$ be the left zero semigroup with elements from $\mathbb{N}$, and let $S$ be the left group $\mathbb{F}_{\setof{a,b}} \times Z$. $S$ is left fairly amenable but is not right fairly amenable.
\Proof{
On the left: let $\xi$ be any finitely-additve measure on $\mathbb{F}_{\setof{a,b}}$ with $\fn\xi{\mathbb{F}_{\setof{a,b}}} = 1$. $\xi$ is necessarily not invariant. Fix an ultrafilter $U$ over $\mathbb{N}$ and define the finitely-additive measure $\mu$ by setting
\[
    \fn\mu{A} \defeq \lim_U \frac1n \sum_{k=1}^n \fn\xi{A\cap (\mathbb{F}_{\setof{a,b}} \times \setof{k}} \quad\text{for all } A\subseteq S.
\]
\begin{pruf}
\step{usual}{$\mu$ exists, is finitely additive, and $\fn\mu{S} = 1$, as usual.}
\step{leftinvar}{$\mu$ is left fairly invariant.}
    \begin{pruf}
    \pf\ Suppose $(g,m)\in S$ acts injectively on the left of $A\subseteq S$: since $Z$ is left zero, this implies that $(x,m_1), (x,m_2)\in A \imp m_1 = m_2$ for all $x\in\mathbb{F}_2$ and $m_1,m_2\in Z$, and thus $\fn\mu{A} = 0$. Then,
    \begin{align*}
        \fn\mu{(g,n)\cdot A} &= \lim_U \frac1n \sum_{k=1}^n \fn\xi{(g,n)A\cap (\mathbb{F}_{\setof{a,b}} \times \setof{k}} \\
            &\le \lim_U \frac1n \fn\xi{\mathbb{F}_{\setof{a,b}}} \\
            &= 0.
    \end{align*}
    \end{pruf}
\end{pruf}

On the right: assume $S$ is right fairly invariant with measure $\nu$.
\begin{pruf}
\step{contradiction}{A contradiction occurs in a similar manner to the usual proof that $\mathbb{F}_2$ is not amenable.}
    \begin{pruf}
    \pf\ Consider one set of words $F(a)\subset \mathbb{F}_{\setof{a,b}}$, which end with the letter $a$. Then 
    \begin{align*}
        (F(a)\times Z)\cdot (a^{-1}, 1) &= (F(a)a^{-1} \times Z) \\    
            &= S\backslash (F(a^{-1}) \times Z), 
    \end{align*}
    and similarly for $F(b)$; hence
    \begin{align*}
        1 &= \fn\nu{S} \\ 
            &= \fn\nu{(F(a) \dotcup F(a^{-1}) \dotcup F(b) \dotcup F(b^{-1}) \dotcup \setof{1}) \times Z} \\
            &\ge \fn\nu{F(a)\times Z} + \fn\nu{F(a^{-1})\times Z} + \fn\nu{F(b)\times Z} + \fn\nu{F(b^{-1})\times Z} \\
            &= \fn\nu{F(a)a^{-1}\times Z} + \fn\nu{F(a^{-1})\times Z} + \fn\nu{F(b)b^{-1}\times Z} + \fn\nu{F(b^{-1})\times Z} \\
            &= \fn\nu{S} + \fn\nu{S} \\
            &= 2,
    \end{align*}
    contradiction. \qed
    \end{pruf}
\end{pruf}

}
\end{example}

\begin{example}\label{ex:freemonogenicinverse}
The free inverse semigroup on one generator $FIS_1$ is fairly amenable both ways. 
\Proof{
From Munn's Theorem on the structure of free inverse semigroups \citep{Lawson1998}, elements of $FIS_1$ can be thought of as triples of integers
\[
FIS_1 \cong \setof{(p,q,r)\in\mathbb{Z}^3: p\ge 0, p+q\ge 0, q+r\ge 0, r\ge 0, p+q+r\ge 0} 
\]
with the product defined by 
\[
(p,q,r)(p',q',r') := \of{\max\setof{p,p'-q}, q+q',\max\setof{r', r-q'}}
\]
for all $(p,q,r), (p',q',r')\in FIS_1$ \citep[p193]{Lawson1998}. 
Consider the increasing sequence given by
\[
    F_n = \setof{(x,y,z)\in FIS_1: x, y, z\le n}.
\]
\begin{pruf}
\step{fnsize}{The sequence $\setof{\abs{F_n}}_{n\in\mathbb{N}}$ is the sequence of ``house numbers'' \citep{oeisA051662}, given by
\[
    \abs{F_n} = (n+1)^3 + \frac16 n(n+1)(2n+1) 
\]
and thus $(n\mapsto \abs{F_n}) \in O(n^3)$.
}
\step{translate}{Let $(p,q,r)\in FIS_1$. By definition,
\[
    (p,q,r)F_n = \setof{\of{\max\setof{p,x-q}, q+y,\max\setof{z, r-y}}: (x,y,z)\in F_n}.
\]}
\step{embiggen}{For large $n$,
\[
    \abs{(p,q,r)F_n} \approx \abs{\setof{\of{x-q, q+y,z}: (x,y,z)\in F_n}}
\]
i.e. the left action of $(p,q,r)$ on $F_n$ is an almost-translation in $\mathbb{Z}^3$, and in particular
\begin{align*}
    \abs{F_n\,\triangle\,(p,q,r)F_n} &\approx \abs{F_n\,\triangle\,\setof{\of{x-q, q+y,z}: (x,y,z)\in F_n}} \\
        &\approx 2qn^2.
\end{align*}
Thus $(n\mapsto \abs{F_n\,\triangle\,(p,q,r)F_n}) \in O(n^2)$,
and therefore the sequence $\setof{F_n}_{n\in\mathbb{N}}$ is \Folner. The right case holds similarly. \qed
}
\end{pruf}
}
\end{example}

Some of the examples and results from above are summarised in Table \ref{table:fairamenabilitysummary}. The variety of interesting examples demonstrate that the ``fair'' modification of invariant finitely-additive measures interacts well with the structure of semigroups. Some important results from group amenability theory are preserved, and examples of fairly amenable semigroups, especially with zeroes, are more gratifying. The given examples of non-fairly amenable semigroups have a certain self-similarity which might be used to create Banach-Tarski-style paradoxes.

\begin{table}
\begin{center}
\begin{tabular}{@{}rll@{}}
\toprule
Kind of semigroup & Classically amenable & Fairly amenable \\
\midrule
Finite & $\iff$ Unique min. ideals & Yes (\ref{cor:allfinite}) \\
With zero & Yes & Sometimes (\ref{cor:addingazero})\\
Monogenic & Yes & Yes (\ref{prop:freecommutativesgp}) \\
Free ($\ge 2$ gen.) & No & No (\ref{ex:freesemigroup2gens}) \\
Abelian & Yes & ? \\
Clifford & Sometimes & Sometimes (\ref{ex:cliffordF2andF1}) \\
Left/right zero sgp & Sided & Yes (\ref{ex:zerosemigroups}) \\
Left/right group & ? & Sometimes (Sided; \ref{ex:leftbutnotright}) \\
Baer-Levi & ? & No (\ref{ex:baerlevi}) \\
Inverse & $\iff$ Max grp hom. im. is & Sometimes \\
Bicyclic & Yes & Yes (\ref{ex:bicyclicmonoid}) \\
Polycyclic & Yes ($\because$ zero) & No (\ref{ex:polycyclicmonoid}) \\
Free monogenic inverse & Yes & Yes (\ref{ex:freemonogenicinverse}) \\
\bottomrule
\end{tabular}
\caption{\label{table:fairamenabilitysummary}Amenability versus fair amenability on different semigroups.}
\end{center}
\end{table}


\section{The convolution partial action}

For real- or complex-valued functions $f:S\mapsto\mathbb{K}$ let the support of $f$ be denoted $\supp(f)$, i.e.
\[
    \supp(f) \defeq \setof{x\in S: f(x) \ne 0}.
\]
When two functions $f$ and $g$ have disjoint support (i.e. $\supp(f) \cap \supp(g) = \emptyset$), we will simply say $f$ and $g$ are \emph{disjoint}.

\begin{figure}
\begin{center}
\begin{tikzpicture}[scale=0.7]
\draw[->] (-3, 0) -- (3, 0) node[anchor=north] {$x\in S$};
\draw[->] (0, -0.5) -- (0, 4);
\draw (-3, 1) -- (-2, 1.3) -- (-0.5, 1) -- (0.5, 3) -- (1,1) -- (2,0.5) -- (2.5,1) -- (3,1) node[anchor=south east] {$f(x)$};
\draw[red,dashed] (-3, 2) -- (3, 2);
\draw[red] (0,2) circle (0.2) node[anchor=south east] {$(0,f(0))$};
\end{tikzpicture}
\begin{tikzpicture}[scale=0.7]
\draw[->] (-3, 0) -- (3, 0) node[anchor=north] {$x\in S$};
\draw[->] (0, -0.5) -- (0, 4);
\draw[red] (-3, 2) -- (3, 2) node[anchor=south east] {$\fn{f}{0x}$};
\end{tikzpicture}
\caption{\label{fig:flattening}The result of the dual left action of 0 on some $f\in\ell^\infty(S)$.}
\end{center}
\end{figure}
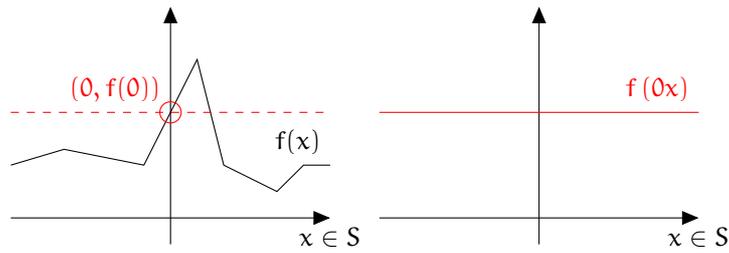

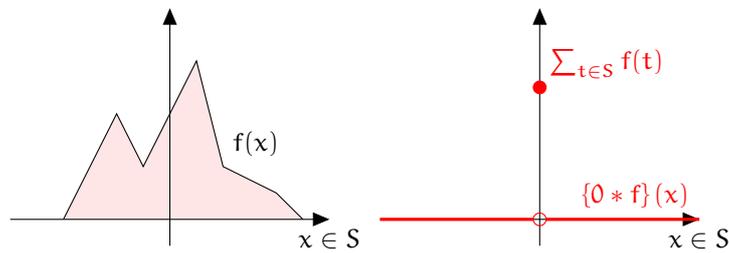
\begin{figure}
\begin{center}
\begin{tikzpicture}[scale=0.7]
\filldraw[fill=red!10] (-2, 0) -- (-1, 2) -- (-0.5, 1) -- (0.5, 3) -- (1,1)  node[anchor=south west] {$f(x)$} -- (2,0.5) -- (2.5,0);
\draw[->] (-3, 0) -- (3, 0) node[anchor=north] {$x\in S$};
\draw[->] (0, -0.5) -- (0, 4);
\end{tikzpicture}
\begin{tikzpicture}[scale=0.7]
\draw[->] (-3, 0) -- (3, 0) node[anchor=north] {$x\in S$};
\draw[->] (0, -0.5) -- (0, 4);
\draw[red,very thick] (-3, 0) -- (-0.12, 0);
\draw[red] (0,0) circle(0.12);
\filldraw[red] (0,2.5) circle (0.12) node[anchor=south west] {$\sum_{t\in S} f(t)$};
\draw[red,very thick] (0.12, 0) -- (3, 0) node[anchor=south east] {$\fn{\setof{0\ast f}}{x}$};
\end{tikzpicture}
\caption{\label{fig:bunchening}The result of the left $\ast$-action of 0 on some $f\in\ell^1(S)$.}
\end{center}
\end{figure}

Recall that convolution of two functions $f,g\in\ell^1(S)$, denoted $f\ast g$, is defined by setting
\[
    \fn{\setof{f\ast g}}{x} \defeq \sum_{st=x} \fn{f}{s} \fn{g}{t} \quad\text{for all }x\in S.
\]
This extends to a left convolution ``action'' of $s\in S$ on $f\in\ell^\infty(S)$, denoted $s\ast f$, which may be defined by setting
\[
	\fn{\setof{s\ast f}}{x} \defeq \sum_{st=x} \fn{f}{t}\quad\text{for all } x\in S.
\]
Alternatively,
\[
    \fn{\setof{s\ast f}}{x} = \sum_{t\in s^{-1}x} \fn{f}{t}\quad\text{for all } x\in S.
\]
For each $s\in S$, let the equivalence relation $\theta_s$ on $S$ be given by setting $x\,\theta_s\,y$ if and only if $sx = sy$, for all $x,y\in S$. Note that each $s^{-1}x$ is precisely a $\theta_s$-equivalence class.

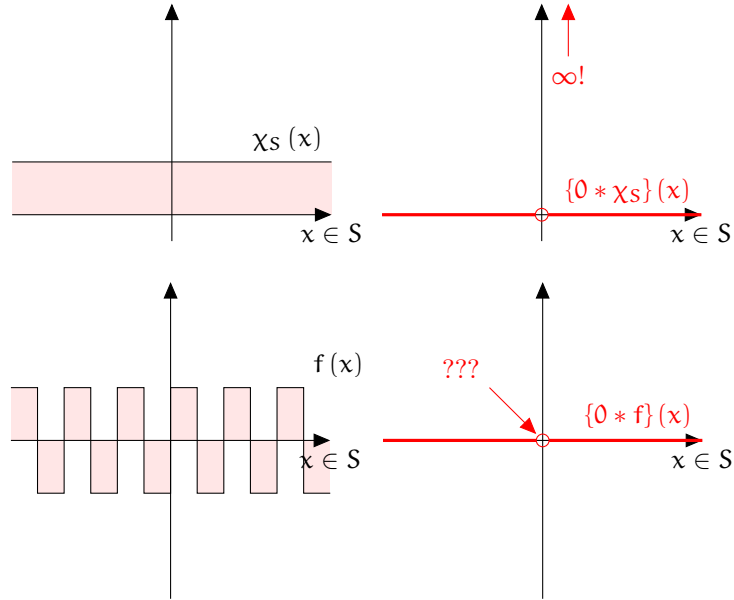
\begin{figure}
    \begin{center}
        \begin{tikzpicture}[scale=0.7]
	        \fill[fill=red!10] (-3, 1) -- (3, 1) -- (3, 0) -- (-3,0) -- (-3, 1);
	        \draw  (-3, 1) -- (3, 1) node[anchor=south east] {$\fn{\chi_S}{x}$};
	        \draw[->] (-3, 0) -- (3, 0) node[anchor=north] {$x\in S$};
	        \draw[->] (0, -0.5) -- (0, 4);
        \end{tikzpicture}
        \begin{tikzpicture}[scale=0.7]
	        \draw[->] (-3, 0) -- (3, 0) node[anchor=north] {$x\in S$};
	        \draw[->] (0, -0.5) -- (0, 4);
	        \draw[red,very thick] (-3, 0) -- (-0.12, 0);
	        \draw[red] (0,0) circle(0.12);
	        \draw[->,red] (0.5, 3) node[anchor=north] {$\infty!$} -- (0.5,4);
	        \draw[red,very thick] (0.12, 0) -- (3, 0) node[anchor=south east] {$\fn{\setof{0\ast \chi_S}}{x}$};
        \end{tikzpicture}
    \end{center}
    \begin{center}
        \begin{tikzpicture}[scale=0.7]
	        \fill[fill=red!10] (-3, 1) -- (-2.5, 1) -- (-2.5, 0) -- (-3,0) -- (-3, 1);
	        \fill[fill=red!10] (-2.5, 0) -- (-2, 0) -- (-2, -1) -- (-2.5,-1) -- (-2.5, 0);
	        \fill[fill=red!10] (-2, 1) -- (-1.5, 1) -- (-1.5, 0) -- (-2,0) -- (-2, 1);
	        \fill[fill=red!10] (-1.5, 0) -- (-1, 0) -- (-1, -1) -- (-1.5,-1) -- (-1.5, 0);
	        \fill[fill=red!10] (-1, 1) -- (-0.5, 1) -- (-0.5, 0) -- (-1,0) -- (-1, 1);
	        \fill[fill=red!10] (-0.5, 0) -- (0, 0) -- (0, -1) -- (-0.5,-1) -- (-0.5, 0);
	        \fill[fill=red!10] (0, 1) -- (0.5, 1) -- (0.5, 0) -- (0,0) -- (0, 1);
	        \fill[fill=red!10] (0.5, 0) -- (1, 0) -- (1, -1) -- (0.5,-1) -- (0.5, 0);
	        \fill[fill=red!10] (1, 1) -- (1.5, 1) -- (1.5, 0) -- (1,0) -- (1, 1);
	        \fill[fill=red!10] (1.5, 0) -- (2, 0) -- (2, -1) -- (1.5,-1) -- (1.5, 0);
	        \fill[fill=red!10] (2, 1) -- (2.5, 1) -- (2.5, 0) -- (2,0) -- (2, 1);
	        \fill[fill=red!10] (2.5, 0) -- (3, 0) -- (3, -1) -- (2.5,-1) -- (2.5, 0);
	        \draw  (-3, 1) -- (-2.5,1) -- (-2.5,-1) -- (-2,-1) -- (-2,1) -- (-1.5,1) -- (-1.5,-1) -- (-1,-1) -- (-1,1) -- (-0.5,1) -- (-0.5,-1) -- (0,-1) -- (0,1) -- (0.5,1) -- (0.5,-1) -- (1,-1) -- (1,1) -- (1.5,1) -- (1.5,-1) -- (2,-1) -- (2,1) -- (2.5,1)  node[anchor=south west] {$\fn{f}{x}$} -- (2.5,-1) -- (3,-1);
	        \draw[->] (-3, 0) -- (3, 0) node[anchor=north] {$x\in S$};
	        \draw[->] (0, -3) -- (0, 3);
        \end{tikzpicture}
        \begin{tikzpicture}[scale=0.7]
	        \draw[->] (-3, 0) -- (3, 0) node[anchor=north] {$x\in S$};
	        \draw[->] (0, -3) -- (0, 3);
	        \draw[red,very thick] (-3, 0) -- (-0.12, 0);
	        \draw[red] (0,0) circle(0.12);
	        \draw[->,red] (-1, 1) node[anchor=south east] {???} -- (-0.12,0.12);
	        \draw[red,very thick] (0.12, 0) -- (3, 0) node[anchor=south east] {$\fn{\setof{0\ast f}}{x}$};
        \end{tikzpicture}
    \caption{\label{fig:convolutionsucks}Some example cases where the convolution partial action of $0$ is not well-defined on $\ell^\infty(S)$.}
    \end{center}
\end{figure}

Unsurprisingly, $\ast$ often fails to be an operation that is closed in $\ell^\infty(S)$, or even well-defined. In contrast to the dual action which ``flattens'' along sections of the domain (see Figure \ref{fig:flattening}), the convolution ``action'' has the appearance of ``bunching up'' the values along the domain (Figure \ref{fig:bunchening}). For an extreme example, suppose $S$ is an infinite semigroup with zero. Then
\[
	0 \ast\chi_S = \sum_{t\in S} \chi_{0\setof{t}} = \sum_{t\in S} \chi_{\setof{0}} = \delta_0,
\]
which takes the ``value'' $\abs{S}=\infty$ at 0. Less extreme cases can also fail to be defined along the entire domain $S$. Examples are depicted in Figure \ref{fig:convolutionsucks}.
There are a few ways this situation might be treated.
\begin{enumerate}[(i)]
\item
We could include, into the scope of discussion, unbounded functions and functions that possibly take the value $\infty$. This makes the $\ast$-``action'' a mapping $S\times l^{\infty}(S) \to \mathbb{C}_{\infty}^S$. This approach is inclusive of degenerate cases such as $\delta_0$, but merely pushes problems relating to singularities into a more complicated place. Additionally this approach does not address those $s\ast f$ which fail to be well-defined, but could still be argued to be bounded.
\item
We could regard convolution as inducing a partial action---simply accept that there will be cases where it is ill-defined, and keep to regions where it is well-defined. Since we wish to apply it to means in $\ell^\infty(S)^*$, we must also keep to the cases that are bounded. It is conceivable that $s\ast f$ exists in $\mathbb{C}^S$ but is unbounded, for instance, $s$ collapses steadily increasing numbers of elements together, but never infinitely many. It is also conceivable that $s\ast f$ is not well-defined because of a failure to converge, but is arguably bounded. 

Now, $s\ast f$ is well-defined and bounded exactly when $s\ast f\in\ell^\infty(S)$. Since $S$ is associative, with $\ell^\infty(S)$ as a set of objects, $S$ induces a set of arrows $A_S$, where for each $s\in S$ there is an arrow from each $f$ to $s\ast f$ wherever $s\ast f\in\ell^\infty(S)$, so $(\ell^\infty(S), A_S)$ defines a semi-category. If $S$ has an identity, then it is a category.
\end{enumerate}
This last point seems interesting, not least because partial actions on C*-algebras are the subject of current research. For our purpose here, we must ask under what conditions is $s\ast f$ bounded, if not $f\in\ell^1(S)$?

\begin{lemma}\label{lemma:injectiveisbounded}
If $s$ acts injectively on the left on $\supp(f)$, then $s\ast f \in \ell^\infty(S)$. 
\Proof{
By hypothesis, $\fn{\setof{s\ast f}}{t}$ is equal to $f(x)$ for some $x\in S$ ($sx = t$) or zero (no such $x$). This is true for any $t\in \supp(f)$, and thus $s\ast f\in \ell^\infty(S)$. \qed
}
\end{lemma}
In particular, $s\ast f$ exists and is bounded whenever $S$ is left cancellative (e.g. is a group). For a semigroup generally, however, the converse does not hold: there may be $f\in \ell^\infty(S)$ such that $s\ast f\in\ell^\infty(S)$ but $s$ is not injective on the support. For example, $f\in\ell^1(\mathbb{N}^0)$ given by $\fn{f}{n} = 2^{-n}$, then $0\ast f = \chi_{\setof{0}}$. 

It is possible to be far more precise than Lemma \ref{lemma:injectiveisbounded} in characterising the elements for which $s\ast f$ is bounded, but it is in essence a restatement of the definition of $\norm{\cdot}_\infty$. For each $s\in S$, those functions $f$ can be thought of as
\begin{enumerate}[(i)]
\item
behaving like elements of $\ell^1(S)$ on the subsets of $S$ on which $s$ does not act injectively, and
\item
 behaving like elements of $\ell^\infty(S)$ on the subsets of $S$ on which $s$ does act injectively.
\end{enumerate}
Since these subsets of $S$ only change with respect to $s$, this suggests, for each $s\in S$, a space\footnote{Not to be confused with either $\ell^p(S)$ or $\ell(S)^*$.} $\ell^{s\ast}(S)$ given by 
\[
    f\in \ell^{s\ast}(S) \iff s\ast f\in\ell^\infty(S).
\]
Clearly, $\ell^1(S) \subseteq \ell^{s\ast}(S) \subseteq \ell^\infty(S)$. 
By definition, $s\ast f\in\ell^\infty(S)$ precisely when there is some fixed finite bound $B \ge \abs{\fn{\setof{s\ast f}}{x}} = \abs{\sum_{t\in s^{-1}x} \fn{f}{t}}$ for all $x\in S$. 

Whether or not $s\ast f\in\ell^\infty(S)$, if $\fn{\setof{s\ast f}}{x}$ is not defined for some $x\in S$, then certainly $x\notin\supp(s\ast f)$, i.e. $\supp(s\ast f)$ can be considered well-defined even if $s\ast f$ is not. Therefore for all $s, f$,  $\supp(s\ast f) = s\cdot\supp(f) \subseteq sS$. Therefore $s\ast\ell^{s\ast}(S)$ can be identified with a subset of $\ell^\infty(sS)$. Since every $f\in\ell^\infty(sS)$, is attainable as some $s\ast g$ for $g\in\ell^{s\ast}(S)$, it follows that
\[
    s\ast\ell^{s\ast}(S) \equiv \ell^\infty(sS);
\]
in particular, $s\ast f\in\ell^\infty(S)$ if, and only if, $(s\ast f)|_{sS} \in \ell^\infty(sS)$. 

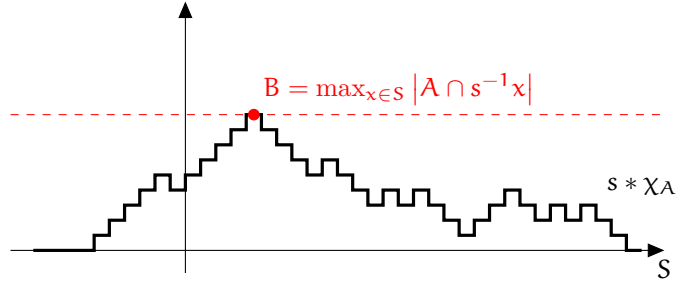
\begin{figure}
\begin{center}
\begin{tikzpicture}
\draw[->] (0,-0.3) -- (0, 3.3);
\draw[->] (-2.3,0) -- (6.3, 0) node[anchor=north] {$S$};
\draw[very thick] (-2.0,0.0) -- (-2.0,0.0) -- (-1.8,0.0) -- (-1.8,0.0) -- (-1.6,0.0) -- (-1.6,0.0) -- (-1.4,0.0) -- (-1.4,0.0) -- (-1.2,0.0) -- (-1.2,0.2) -- (-1.0,0.2) -- (-1.0,0.4) -- (-0.8,0.4) -- (-0.8,0.6) -- (-0.6,0.6) -- (-0.6,0.8) -- (-0.4,0.8) -- (-0.4,1.0) -- (-0.2,1.0) -- (-0.2,0.8) -- (-0.0,0.8) -- (-0.0,1.0) -- (0.2,1.0) -- (0.2,1.2) -- (0.4,1.2) -- (0.4,1.4) -- (0.6,1.4) -- (0.6,1.6) -- (0.8,1.6) -- (0.8,1.8) -- (1.0,1.8) -- (1.0,1.6) -- (1.2,1.6) -- (1.2,1.4) -- (1.4,1.4) -- (1.4,1.2) -- (1.6,1.2) -- (1.6,1.0) -- (1.8,1.0) -- (1.8,1.2) -- (2.0,1.2) -- (2.0,1.0) -- (2.2,1.0) -- (2.2,0.8) -- (2.4,0.8) -- (2.4,0.6) -- (2.6,0.6) -- (2.6,0.8) -- (2.8,0.8) -- (2.8,0.6) -- (3.0,0.6) -- (3.0,0.8) -- (3.2,0.8) -- (3.2,0.6) -- (3.4,0.6) -- (3.4,0.4) -- (3.6,0.4) -- (3.6,0.2) -- (3.8,0.2) -- (3.8,0.4) -- (4.0,0.4) -- (4.0,0.6) -- (4.2,0.6) -- (4.2,0.8) -- (4.4,0.8) -- (4.4,0.6) -- (4.6,0.6) -- (4.6,0.4) -- (4.8,0.4) -- (4.8,0.6) -- (5.0,0.6) -- (5.0,0.4) -- (5.2,0.4) -- (5.2,0.6) -- (5.4,0.6) node[anchor=south west] {$s\ast\chi_A$} -- (5.4,0.4) -- (5.6,0.4) -- (5.6,0.2) -- (5.8,0.2) -- (5.8,0.0) -- (6.0,0.0);
\draw[dashed,red] (-2.3,1.8) -- (6.3,1.8);
\filldraw[red] (0.9, 1.8) circle (0.07) node[anchor=south west] {$B = \max_{x\in S}\abs{A\cap s^{-1}x}$};
\end{tikzpicture}
\caption{\label{figure:sampleactiononindicator}Diagram accompanying Lemma \ref{lemma:mmmmonsterkill}.}
\end{center}
\end{figure}

\begin{lemma}\label{lemma:mmmmonsterkill}
For all $s\in S$ and $A\subseteq S$, the following conditions are equivalent.
\begin{enumerate}[(i)]
\item $s\ast\chi_A\in\ell^\infty(S)$.
\item There exists a finite partition $\setof{A_i}_{i\in I}$ of $A$ such that $s$ acts injectively on the left of each $A_i$.
\item $s\ast\chi_A$ is simple.
\end{enumerate}
\Proof{
\begin{enumerate}
\item[(i) $\imp$ (ii):]
Suppose $\norm{s\ast\chi_A}_\infty = B<\infty$. $B$ is a non-negative integer which $s\ast\chi_A$ attains, since the value at each point is a sum of values in $\setof{0,1}$. For all $x\in S$ we have
\begin{align*}
    \fn{\setof{s\ast\chi_A}}{x} &= \sum_{t\in s^{-1}x} \fn{\chi_A}{t} \\
        &= \abs{A\cap s^{-1}x} \\
        &\le B \quad\text{by hypothesis.}
\end{align*}
For $i=1,\ldots,B$ let $A_i$ consist of one choice element from each $(A\cap s^{-1}x)\backslash \bigcup_{j<i} A_j$ for $x\in S$ (where it is not empty).\footnote{The Axiom of Choice is not required because the set is finite.} Then $B$ choices are made, each $A\cap s^{-1}x$ is exhausted, and $I = \setof{1,\ldots,B}$ is finite. The finite collection $\setof{A_i}_{i\in I}$ is a partition of $A$, since the sets of the form $s^{-1}x$ for each $x\in S$ are either empty, or distinct $\theta_s$-equivalence classes. $s$ acts injectively on the left of each $A_i$, as $A_i\cap s^{-1}x$ is either empty or a singleton set.

\item[(ii) $\imp$ (iii):] Suppose there is a finite partition $\setof{A_i}_{i\in I}$ of $A$ such that $s$ acts injectively on the left of $A_i$. Then $\chi_A = \sum_{i\in I} \chi_{A_i}$ and $s\ast\chi_{A_i} = \chi_{sA_i}$ for each $i\in I$, and thus 
\[
    s\ast\chi_A = s\ast\sum_{i\in I} \chi_{A_i} = \sum_{i\in I} s\ast\chi_{A_i} = \sum_{i\in I} \chi_{sA_i},
\]
which is a linear combination of finitely-many indicator functions, i.e. is simple.

\item[(iii) $\imp$ (i):] If $s\ast\chi_A$ is simple then by definition it consists of a linear combination of finitely-many indicator functions, and thus attains some finite bound. \qed
\end{enumerate}
}
\end{lemma}

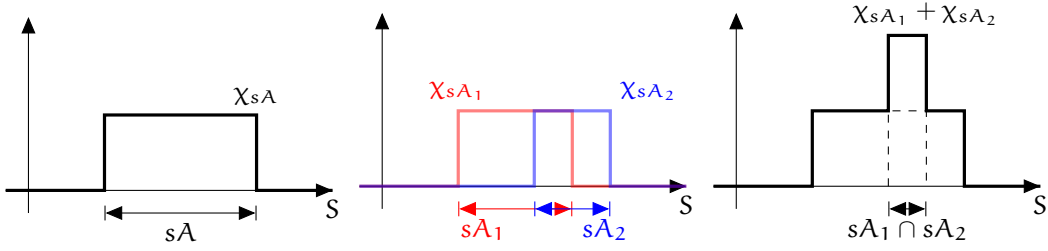
\begin{figure}
\begin{center}
\begin{tikzpicture}
\draw[->] (0,-0.3) -- (0, 2.3);
\draw[->] (-0.3,0) -- (4, 0) node[anchor=north] {$S$};
\draw[very thick] (-0.3,0) -- (1,0) -- (1,1) -- (3,1) node[anchor=south] {$\chi_{sA}$} -- (3,0) -- (4,0);
\draw[|<->|] (1,-0.3) -- node[anchor=north] {$sA$} (3,-0.3);
\end{tikzpicture}
\begin{tikzpicture}
\draw[->] (0,-0.3) -- (0, 2.3);
\draw[->] (-0.3,0) -- (4, 0) node[anchor=north] {$S$};
\draw[very thick,red,draw opacity=0.5] (-0.3,0) -- (1,0) -- (1,1) node[anchor=south] {$\chi_{sA_1}$} -- (2.5,1) -- (2.5,0) -- (4,0);
\draw[very thick,blue,draw opacity=0.5] (-0.3,0) -- (2,0) -- (2,1) -- (3,1) node[anchor=south west] {$\chi_{sA_2}$} -- (3,0) -- (4,0);
\draw[|<->|,red] (1,-0.3) -- node[anchor=north east] {$sA_1$} (2.5,-0.3);
\draw[|<->|,blue] (2,-0.3) -- node[anchor=north west] {$sA_2$} (3,-0.3);
\end{tikzpicture}
\begin{tikzpicture}
\draw[->] (0,-0.3) -- (0, 2.3);
\draw[->] (-0.3,0) -- (4, 0) node[anchor=north] {$S$};
\draw[very thick] (-0.3,0) -- (1,0) -- (1,1) -- (2,1) -- (2,2) -- (2.5,2) node[anchor=south] {$\chi_{sA_1} + \chi_{sA_2}$} -- (2.5,1) -- (3,1) -- (3,0) -- (4,0);
\draw[|<->|] (2,-0.3) -- node[anchor=north] {$sA_1\cap sA_2$} (2.5,-0.3);
\draw[dashed] (2,0) -- (2,1) -- (2.5,1) -- (2.5,0);
\end{tikzpicture}
\caption{\label{fig:chi1chi2}An example of $s\ast\chi_A \ge \chi_{sA}$. $\chi_{sA}\le\chi_{sA_1} + \chi_{sA_2}$, where $A = A_1\dotcup A_2$ and $s$ acts injectively on $A_1$ and $A_2$ but not $A$ as a whole.}
\end{center}
\end{figure}

Another impediment to deducing standard results includes the difficulty in working with even simple functions. Suppose $f\in \ell^\infty(S)$ is simple, and thus there exists a finite index set $I$, and collections of numbers $\setof{a_i\in\mathbb{C}: i\in I}$ and sets $\setof{A_i\in\Sc{P}(S): i\in I}$ such that
\(
	f = \sum_{i\in I} a_i\chi_{A_i}.
\)
Where it exists, $\ast$ distributes over $+$, and clearly if $s\ast f$ is bounded then $s\ast\chi_{A_i}$ is also bounded for each $i\in I$. Therefore,
\[
	s\ast f = \sum_{i\in I} a_i \cdot (s\ast\chi_{A_i}).
\]
However, if the action of $s$ is \emph{not} injective on each $A_i$, it isn't at all likely that
\[
    \sum_{i\in I} a_i\cdot (s\ast\chi_{A_i}) \ge \sum_{i\in I} a_i \chi_{sA_i}
\]
is saturated, and in fact $\sum_{i\in I} a_i \chi_{sA_i}$ could vary depending upon the selection of $\setof{A_i}_{i\in I}$.

Fortunately, if $s\ast f$ is bounded then each $s\ast\chi_{A_i}$ is bounded and therefore by Lemma \ref{lemma:mmmmonsterkill} is simple, and also, there exists a finite partition $\setof{B_{ij}}_{j\in J_i}$ of each $A_i$ such that $s$ acts injectively on the left of each $B_{ij}$. Thus
\[
    s\ast f = \sum_{i\in I} a_i \sum_{j\in J_i} \chi_{sB_{ij}},
\]
and hence if $f$ is simple then so is $s\ast f$ (where it exists and is bounded).

\section{Integrating $s\ast f$}

\begin{definition}\label{definition:astinvariantwherebounded}
Let $m\in\ell^\infty(S)^*$. $m$ is \emph{left $\ast$-invariant} if
\[
    \fn{m}{f} = \fn{m}{s\ast f}  
\]
for all $s\in S$ and $f\in\ell^\infty(S)$ wherever $s\ast f\in\ell^\infty(S)$.
\end{definition}

The purpose of this section will be to show that Definition \ref{definition:astinvariantwherebounded} is equivalent to left fair amenability of $S$, i.e. the existence of a left $\ast$-invariant (where bounded) mean is equivalent to the existence of a left fairly-invariant probability measure.

Suppose $S$ supports a left $\ast$-invariant mean $m$ as described in Definition \ref{definition:astinvariantwherebounded}. It is easy to see why Definition \ref{definition:astinvariantwherebounded} is at least as strong as left fair amenability: when $s$ acts injectively on the left of $A$, $s\ast\chi_A = \chi_{sA}$, and so a $\ast$-invariant mean can be applied to the indicator functions. To show the converse, I shall integrate with respect to $\mu$.

First, let us consider indicator functions.

\begin{lemma}\label{lemma:indicatorastinvariant}
Let $S$ be a left fairly amenable semigroup with measure $\mu$, $s\in S$, and $A\subseteq S$. If $s\ast\chi_A\in\ell^\infty(S)$ then 
\[
    \int\of{s\ast\chi_A}\D\mu = \int\chi_A\D\mu.
\]
\Proof{
If $s\ast\chi_A\in\ell^\infty(S)$, then there is the finite partition $\setof{A_i}_{i\in I}$ of $A$ provided by Lemma \ref{lemma:mmmmonsterkill} such that $s$ acts injectively on the left of each $A_i$. Then
\begin{align*}
    \int\of{s\ast\chi_A}\D\mu &= \int\of{\sum_{i\in I}\chi_{sA_i}}\D\mu \\
        &= \sum_{i\in I}\fn\mu{sA_i} \quad\text{by definition} \\
        &= \sum_{i\in I}\fn\mu{A_i} \quad\because\text{fair invariance} \\
        &= \fn\mu{A} \quad\because\text{finitely additive}\\
        &= \int\chi_A\D\mu \quad\text{again by definition,}
\end{align*}
as required.\qed}
\end{lemma}

\begin{lemma}\label{lemma:simpleastinvariant}
Let $S$ be a left fairly amenable semigroup with measure $\mu$, $s\in S$, and $f\in\ell^\infty_+(S)$ is a simple function. If $s\ast f\in\ell^\infty(S)$ then 
\[
    \int\of{s\ast f}\D\mu = \int f\D\mu.
\]
\Proof{
There are the requisite finite index set $I$, sets $\setof{A_i}_{i\in I}$ and values $a_i\in\mathbb{R}^+$ for $i\in I$ such that
\(
	f = \sum_{i\in I} a_i\chi_{A_i}. 
\)
If $s\ast f\in\ell^\infty(S)$ then $s\ast\chi_{A_i}\in\ell^\infty(S)$ for each $i\in I$, and therefore
\begin{align*}
    \int\of{s\ast f}\D\mu &= \int\of{s\ast\sum_{i\in I} a_i \chi_{A_i}}\D\mu\\
        &= \int\of{\sum_{i\in I} a_i\cdot\of{s\ast\chi_{A_i}}}\D\mu\\
        &= \sum_{i\in I} a_i \int\of{s\ast\chi_{A_i}}\D\mu \\
        &= \sum_{i\in I} a_i \int\chi_{A_i}\D\mu \quad\because\text{Lemma \ref{lemma:indicatorastinvariant}}\\
        &= \int\of{\sum_{i\in I} a_i\chi_{A_i}}\D\mu \\
        &= \int f\D\mu
\end{align*}
as required.\qed}
\end{lemma}


Let $\ell^\infty_+(S)$ denote the subset of $\ell^\infty(S)$ consisting of bounded real-valued non-negative functions on $S$. 

Not every simple function $h\le s\ast f$ is of the form $s\ast g$ for a simple $g\in\ell^\infty_+(S), g\le f$. For example, let $f\in\ell^1(\mathbb{N}^0)$ with $f(n) = 1/2^n$ for each $n$. Then $\setof{0\ast f}(0) = 2$. $h=0\ast f$ itself is simple. However, there is no $g\le f$ such that $g$ is simple and $0\ast g = h$ (that would require either $g$ to be non-simple or $g>f$). It is nevertheless sufficient to range over all functions of the form $s\ast g$ for simple $g\le f$, when integrating $s\ast f$, as $s\ast g$ approximates $s\ast f$ increasingly well as $g$ gains detail.

\begin{lemma}\label{lemma:positiveastinvariant}
Let $S$ be a left fairly amenable semigroup with measure $\mu$, $s\in S$, and $f\in\ell^\infty_+(S)$. If $s\ast f\in\ell^\infty_+(S)$ then 
\[
    \int\of{s\ast f}\D\mu = \int f\D\mu.
\]
\Proof{
If $s\ast f\in\ell^\infty_+(S)$ then $s\ast h\in\ell^\infty_+(S)$ for every simple function $h\le f$. Thus
\begin{align*}
    \int\of{s\ast f}\D\mu &= \sup\setof{\int h\D\mu: h\le (s\ast f), h\text{ is simple}} \quad\text{by definition}\\
        &= \sup\setof{\int (s\ast g)\D\mu: (s\ast g)\le (s\ast f), g\text{ is simple}}\\
        &= \sup\setof{\int g\D\mu: g\le f, g\text{ is simple}}\quad\because\text{Lemma \ref{lemma:simpleastinvariant}}\\
        &= \int f\D\mu \quad\text{by definition}
\end{align*}
as required.\qed
}
\end{lemma}
The next lemma is is entirely routine.
\begin{lemma}\label{lemma:arbitraryastinvariant}
Let $S$ be a left fairly amenable semigroup with measure $\mu$, $s\in S$, and real-valued $f\in\ell^\infty(S)$. If $s\ast f\in\ell^\infty(S)$ then 
\[
    \int\of{s\ast f}\D\mu = \int f\D\mu.
\]
\Proof{
There exist $f^+, f^- \in\ell^\infty_+(S)$ such that $f = f^+ - f^-$. If $s\ast f\in\ell^\infty(S)$ then so too $s\ast f^+$ and $s\ast f^-$, thus
\begin{align*}
    \int (s\ast f)\D\mu &= \int (s\ast (f^+ - f^-))\D\mu \\
        &= \int ((s\ast f^+) - (s\ast f^-))\D\mu \\
        &= \int (s\ast f^+) \D\mu - \int (s\ast f^-)\D\mu \\
        &= \int f^+\D\mu - \int f^-\D\mu \quad\because\text{Lemma \ref{lemma:positiveastinvariant}}\\
        &= \int (f^+ - f^-)\D\mu \\
        &= \int f\D\mu
\end{align*}
as required. \qed
}
\end{lemma}
The complex-valued case is even more pedestrian, so it is omitted.

\begin{theorem}[Main Theorem]
A semigroup $S$ is left fairly amenable if, and only if, there exists a left $\ast$-invariant mean in $\ell^\infty(S)^*$.
\Proof{
Suppose $S$ is left fairly amenable with the finitely-additive measure $\mu$. By Lemma \ref{lemma:arbitraryastinvariant}, the integral with respect to $\mu$ is $\ast$-invariant, therefore use the mean $m\in\ell^\infty(S)^*$ given by setting
\[
    \fn{m}{f} \defeq \int f\D\mu \quad\text{for all } f\in\ell^\infty(S).
\]
Conversely, if $S$ supports a left $\ast$-invariant mean $m$, define the measure $\mu\in[0,1]^{\Sc{P}(S)}$ by setting 
\[
    \fn\mu{A} \defeq \fn{m}{\chi_A}\quad\text{for all } A\in\Sc{P}(S). 
\]
Then, if $s\in S$ acts injectively on the left of $A\in\Sc{P}(S)$, $s\ast\chi_A = \chi_{sA}$, and then
\begin{align*}
    \fn\mu{sA} &= \fn{m}{\chi_{sA}} \\
        &= \fn{m}{s\ast\chi_A} \\
        &= \fn{m}{\chi_A} = \fn\mu{A},
\end{align*}
as required. \qed}
\end{theorem}


\subsection{Acknowledgements}
I extend a heartfelt thank-you to my supervisor Dr Des FitzGerald 
 for his continued patience and wisdom during the project.


\bibliographystyle{kluwer}
\bibliography{library,OEIS}

\end{document}